\documentclass[reqno]{amsart}
\usepackage{latexsym,amsmath,amsfonts,amscd,amssymb}

\usepackage{eepic}
\theoremstyle{plain}
\newtheorem{theorem}{Theorem}[section]
\newtheorem{corollary}[theorem]{Corollary}
\newtheorem{lemma}[theorem]{Lemma}
\newtheorem{proposition}[theorem]{Proposition}
\theoremstyle{definition}
\newtheorem{definition}[theorem]{Definition}
\newtheorem{remark}[theorem]{Remark}
\numberwithin{equation}{section}

\renewcommand{\a}{\alpha}

\newcommand{\s}{\sigma}

\newcommand{\cN}{\mathcal{N}}

\newcommand{\cO}{\mathcal{O}}

\newcommand{\cS}{\mathcal{S}}

\newcommand{\RR}{\mathbb{R}}
\newcommand{\CC}{\mathbb{C}}
\newcommand{\HH}{\mathbb{H}}
\newcommand{\PP}{\mathbb{P}}
\newcommand{\ZZ}{\mathbb{Z}}

\newcommand{\inc}{\hookrightarrow}
\newcommand{\lto}{\longrightarrow}

\newcommand{\x}{\times}
\newcommand{\ox}{\otimes}
\newcommand{\iso}{\cong}

\newcommand{\coeff}{\mathop{\mathrm{coeff}}}
\newcommand{\Sym}{\mathrm{Sym}}
\newcommand{\Res}{\mathrm{Res}}
\newcommand{\GCD}{\mathrm{gcd}}
\DeclareMathOperator{\Jac}{Jac}
\DeclareMathOperator{\rk}{rk} 

\DeclareMathOperator{\Hom}{Hom} \DeclareMathOperator{\Ext}{Ext}
 \DeclareMathOperator{\Gr}{Gr}
\DeclareMathOperator{\GL}{GL}\DeclareMathOperator{\Aut}{Aut}

\newcommand{\scp}{{\s_c^+}}
\newcommand{\smp}{{\s_m^+}}
\newcommand{\scm}{{\s_c^-}}
\newcommand{\scpm}{{\s_c^\pm}}
\newcommand{\moduli}{\cN_\s}

\newcommand{\modulimas}{\cN_{\s_c^+}}
\newcommand{\modulimenos}{\cN_{\s_c^-}}
\newcommand{\Smas}{\cS_{\scp}}
\newcommand{\Smenos}{\cS_{\scm}}

\newcommand{\modulimmas}{\cN_{\smp}}

\title[Hodge polynomials of moduli spaces of triples]{Hodge polynomials
of the moduli spaces \\ of triples of rank~$(2,2)$}


\author{Vicente Mu\~noz}
  \address{Departamento de Matem\'aticas \\
  CSIC \\ Serrano 113 bis
  \\ 28006 Madrid \\ Spain \vspace{-3mm}}
  \address{Institute for Advanced Study \\ Einstein Drive \\
  Princeton, New Jersey 08540 \\ USA}
  \email{vicente.munoz@imaff.cfmac.csic.es}

\author{Daniel Ortega}
  \address{Departamento de Matem\'aticas \\
  Universidad Aut\'onoma de Madrid
  \\ 28049 Madrid \\ Spain}
  \email{daniel.ortega@uam.es}

\author{Maria-Jes{\'u}s V{\'a}zquez-Gallo}
  \address{Departamento de Ingenier\'{\i}a Civil: Servicios
  Urbanos\\
    Unidad Docente: Matem\'aticas \\
    Escuela de Ingenier{\'\i}a de Obras P{\'u}blicas \\
    Universidad Polit{\'e}cnica de Madrid \\
    Alfonso XII 3 y 5 \\ 28014 Madrid \\ Spain}
  \email{mariajesus.vazquez@upm.es}

\begin{document}

\begin{abstract}
 Let $X$ be a smooth projective curve of genus $g\geq 2$ over the
 complex numbers. A holomorphic triple $(E_1,E_2,\phi)$ on $X$ consists of two
 holomorphic vector bundles $E_{1}$ and $E_{2}$ over $X$ and a
 holomorphic map $\phi \colon E_{2} \to E_{1}$.
 There is a concept of stability for triples which depends on a
 real parameter $\sigma$.
 In this paper, we determine the Hodge polynomials of the moduli
 spaces of $\sigma$-stable triples with $\rk(E_1)=\rk(E_2)=2$, using the theory of mixed Hodge
 structures (in the cases that they are smooth and compact).
 This gives in particular the Poincar{\'e} polynomials
 of these moduli spaces.
 As a byproduct, we also give the Hodge
 polynomial of the moduli space of even degree rank $2$ stable
 vector bundles.
\end{abstract}

\maketitle

\section{Introduction}
\label{sec:introduction}

Let $X$ be a smooth projective curve of genus $g\geq  2$ over the
field of complex numbers. A holomorphic triple $T =
(E_{1},E_{2},\phi)$ on $X$ of rank $(n_1,n_2)$ consists of two
holomorphic vector bundles $E_{1}$ and $E_{2}$ over $X$ (of ranks
$n_1$ and $n_2$, and degrees $d_1$ and $d_2$, respectively) and a
holomorphic map $\phi \colon E_{2} \to E_{1}$. There is a concept
of stability for a triple which depends on the choice of a
parameter $\s \in \RR$. Let $\moduli$ and $\moduli^s$ denote the
moduli spaces of $\sigma$-semistable and $\sigma$-stable
triples, respectively. These have been widely studied in
\cite{BGP,BGPG,GPGM,MOV}.

The range of the parameter $\s$ is an interval $I\subset \RR$
split by a finite number of \textit{critical values} $\s_c$ in
such a way that, when $\s$ moves without crossing a critical
value, then $\moduli$ remains unchanged, but when $\s$ crosses a
critical value, $\moduli$ undergoes a transformation which we call
\textit{flip}. The study of this process allows us to obtain
information on the topology of all moduli spaces $\moduli$, for
any $\s$, once we know such information for one particular
$\moduli$ (usually the one corresponding to the minimum or maximum
possible value of the parameter).

One of the main motivations to study the topology of the moduli
spaces of triples is that they appear when looking at the topology
of the moduli spaces of Higgs bundles \cite{Hi,Go,GPGM} via Morse
theory techniques. Higgs bundles are pairs $(E,\Phi)$, formed by a
holomorphic vector bundle $E$ of rank $r$ and a holomorphic map
$\Phi:E\to E\otimes K$, where $K$ is the canonical bundle of the
curve, and they are intimately related to the representation
varieties of the fundamental group of the surface underlying the
complex curve into the general Lie group $\GL(r,\CC)$. The moduli
spaces of triples, and the more general moduli spaces of chains
\cite{AG,GK,AGS} appear as critical sets of a natural Morse-Bott
function on the moduli space of Higgs bundles \cite{Hi,Go}.

When the rank of $E_2$ is one, we have the so-called
\textit{pairs} \cite{BD,GP,MOV}. The moduli spaces of pairs are
smooth for any rank $n_1$, and in the case of rank $n_1=2$ and
fixed determinant, they are very well-understood thanks to the
work of Thaddeus \cite{Th}. In this case, the flips have a very
nice geometrical interpretation, consisting of blowing up an
embedded subvariety and then blowing-down the exceptional divisor
in a different way. Moreover, there are also very explicit
descriptions of the moduli spaces of pairs for the minimum and
maximum possible values of $\s$.

The flips do not have such a nice behavior for moduli spaces of triples of rank
$(n_1,n_2)$ with \mbox{$n_1+n_2>3$}. The flip locus may have singularities, it may
consist of several irreducible components intersecting in a non-transverse way, the
moduli spaces themselves may have singularities for \mbox{$n_1,n_2\geq 2$}, and the
moduli spaces for $\s$ large are difficult to handle in the situation when
$n_1=n_2$, since then they are described in terms of Quot schemes.

These difficulties can be overcome in two different ways. The
first way is to introduce parabolic structures with generic
weights. The moduli spaces of parabolic triples have been studied
in \cite{GPGM}, where the Poincar{\'e} polynomials have been given for
the moduli of parabolic triples of ranks $(2,1)$. The parabolic
weights tend to prevent the singularities of the moduli spaces and
flip loci. However, for obtaining information on the moduli space
of non-parabolic triples, one should relate the parabolic and the
non-parabolic situations.

The second route to compute the Poincar{\'e} polynomials of the moduli spaces of
triples was introduced in \cite{MOV}. It consists of using the theory of mixed Hodge
structures of Deligne \cite{De} to compute the Hodge polynomials of the moduli
space. The Hodge polynomials recover the usual Poincar{\'e} polynomial when we deal
with a smooth compact algebraic variety, but they can be defined for non-smooth and
non-compact algebraic varieties as well. This allows to compute the Poincar{\'e}
polynomials of the moduli spaces of triples which are smooth and compact, no matter
if the flip loci have singularities.

In this paper, we use mixed Hodge theory to compute the Hodge
polynomials of some of the moduli spaces of triples of rank $(2,2)$.
By the results of \cite{BGPG}, if $d_1-d_2 > 4g-4$ then
$\cN_\s^s$ is smooth. Moreover, when $d_1+d_2$ is odd,
the moduli spaces $\cN_\s$ only consist of $\s$-stable
triples for non-critical values of $\sigma$, therefore
$\cN_\s^s$ are projective varieties. Because of this, we shall
compute the Hodge polynomials of the moduli spaces of triples
of rank $(2,2)$ in the case $d_1-d_2 > 4g-4$ and $d_1+d_2$ odd.
This gives in particular the Poincar{\'e} polynomials of these
moduli spaces.

We start by reviewing the rudiments of mixed Hodge theory and the standard results
on triples that we shall use throughout the paper, in Sections \ref{sec:virtual} and
\ref{sec:stable-triples}. Then Section \ref{sec:Poincare(2,1)} recalls the
computations of the Hodge polynomials of the moduli spaces of triples of ranks
$(2,1)$ and $(1,2)$, from \cite{MOV}. In Section~\ref{sec:poly-bundles-rk-two} we
use the Hodge polynomial of the moduli spaces of triples to deduce the Hodge
polynomials of the moduli spaces of rank $2$ stable vector bundles. The case of odd
degree rank $2$ bundles is already known \cite{BaR,EK,MOV}, but we do the case of
even degree rank $2$ stable bundles, proving the following result (see~Theorem
\ref{thm:rank2even}).

\noindent \textbf{Theorem A.}\ {\em Let $M^s(2,d)$ denote the moduli
space of rank $2$, degree $d$ stable vector bundles on $X$. If $d$ is even
then the Hodge polynomial of $M^s(2,d)$ is
 \begin{align*}
   e(M^s(2,d))=\, &
  \frac{1}{2(1-uv)(1-(uv)^2)} \bigg(
  2(1+u)^{g}(1+v)^g(1+u^2v)^{g}(1+uv^2)^g   \\ & -
  (1+u)^{2g}(1+v)^{2g} (1+ 2 u^{g+1} v^{g+1} -u^2v^2) -
  (1-u^2)^g(1-v^2)^g(1-uv)^2 \bigg) \, .
 \end{align*}
}
 \medskip

Note that the moduli space $M^s(2,d)$ is smooth but non-compact.

Next we move to the study of the moduli spaces of
triples of rank $(2,2)$, which are the main focus of the paper. The
critical values are computed in Section \ref{sec:critical(2,2)}.

In Section \ref{sec:small} we compute the Hodge
polynomial of the moduli space of stable triples of rank $(2,2)$ for
the smallest allowable values of the parameter $\s$, proving the following
result (see Theorem \ref{thm:(2,2,odd,even)} and Corollary
\ref{cor:(2,2,odd,even)-2}).

\noindent \textbf{Theorem B.}\ 
 {\em Let $\cN_\s=\cN_\s(2,2,d_1,d_2)$ be the moduli space
  of $\s$-stable triples of rank $(2,2)$. Assume that $d_1-d_2 > 4g-4$ and
  $d_1+d_2$ is odd. Let $\sigma_m=\frac{d_1}2 - \frac{d_2}2$ be the minimum
  value of the parameter $\sigma$ and $\smp=\sigma_m + \epsilon$ for
  $\epsilon>0$ small. Then  $\cN_\smp$ is smooth and projective, it only
  consists of stable triples, and its Hodge polynomial is
  $$
 \begin{aligned}
 e(\modulimmas) = & \ \frac{(1+u)^{2g}(1+v)^{2g}(1-(uv)^N)
 (u^gv^g(1+u)^g(1+v)^g-(1+u^2v)^g(1+uv^2)^g)}{(1-uv)^3(1-(uv)^2)^2}
 \cdot \\ & \bigg(
  (1+u)^{g} (1+v)^g (u^{g+1}v^{g+1} + u^{N+g-1}v^{N+g-1}) -
  (1+u^2v)^{g}(1+uv^2)^g (1 + u^Nv^{N})\bigg) \, ,
 \end{aligned}
  $$
where $N=d_1-d_2-2g+2$. }

\smallskip

\enlargethispage{1\baselineskip}

Under the condition $d_1-d_2>4g-4$, the Hodge polynomial of $\cN_\smp(2,2,d_1,d_2)$
when both $d_1,d_2$ are odd is easily given (see Theorem \ref{thm:(2,2,odd,odd)}).
When both $d_1,d_2$ are even, it may be computed with similar techniques to those of
Theorem \ref{thm:(2,2,odd,even)}. However, to remove the condition $d_1-d_2>4g-4$ is
not possible with the current techniques.

The contribution of the flips to the Hodge polynomials of the moduli spaces of
$\s$-stables triples of rank $(2,2)$ is computed in Section
\ref{sec:simple}. This is added up to the information for the Hodge polynomial
of the small parameter moduli space to get the Hodge polynomial of the
moduli space of $\s$-stable triples of rank $(2,2)$ for the largest
values of $\s$ in Section \ref{sec:large}. We get the following result
(see Theorem \ref{thm:finally} and Corollary \ref{cor:finally}).

\noindent \textbf{Theorem C.}\ 
{\em  Let $\cN_\s=\cN_\s(2,2,d_1,d_2)$ be the moduli space
  of $\s$-stable triples of rank $(2,2)$. Assume that $d_1-d_2 > 4g-4$ and
  $d_1+d_2$ is odd. Let $\sigma_M=d_1 - d_2$. Then all the moduli spaces
  $\cN_\s$ are isomorphic for $\s>\s_M$. Let
  $\s_M^+=\sigma_M + \epsilon$ for
  $\epsilon>0$. Then
  $\cN_{\s_M^+}$ is smooth and projective, it only
  consists of stable triples, and
 $$
 \begin{aligned}
 e(\cN_{\s_M^+})=& \ \frac{(1+u)^{2g}(1+v)^{2g}}{(1-uv)^3(1-(uv)^2)^2}
    \Bigg[
    (1+u^2v)^{2g}(1+uv^2)^{2g}(1-(uv)^{2N})\\
    &
    -N \, (1+u^2v)^g(1+uv^2)^g(1+u)^g(1+v)^g (uv)^{N+g-1}(1-(uv)^2)
    \\
    &
    +
    (1+u)^{2g}(1+v)^{2g}(1+uv)^2(uv)^{2g -2 +(N+1)/2 } 
    \Big((1-(uv)^{N+1}) - \frac{N+1}{2} \, (1-uv)(1+(uv)^{N})\Big)\\
    &
    -g(1+u)^{2g-1}(1+v)^{2g-1}
    (1-(uv)^2)^2(uv)^{2g -2 +(N+1)/2}(1-(uv)^{N})
    \Bigg] \, ,
\end{aligned}
$$
where $N=d_1-d_2-2g+2$.
}
\medskip

The computation of the contribution of the flips to the Hodge polynomials of the
moduli spaces of $\s$-stables triples of rank $(2,2)$ is done under the assumptions
$d_1+d_2$ odd and $d_1-d_2>2g-2$. This~can be extended to the case $d_1+d_2$ even,
keeping in mind that in this case we will find the Hodge polynomials of the moduli
spaces $\cN_\s^s$ which are non-compact and of the moduli spaces $\cN_\s$ which have
singularities at non-stable points. However the assumption $d_1-d_2>2g-2$ cannot be
removed with the current techniques.

The Poincar\'e polynomials of the moduli spaces $\cN_{\s_m^+}$ and $\cN_{\s_M^+}$
are obtained from the Hodge polynomials, for $d_1-d_2 > 4g-4$ and
$d_1+d_2$ odd (see Corollaries \ref{cor:(2,2,odd,even)} and
\ref{cor:finally-n}), since they are smooth projective varieties.

{\bf Acknowledgements:} First and third authors partially supported
through grant MCyT (Spain) MTM2004-07090-C03-01/02. 

\section{Hodge Polynomials}
\label{sec:virtual}

\subsection{Hodge-Deligne theory}
\label{subsec:Hodge-Deligne}

Let us start by recalling the Hodge-Deligne theory of algebraic
varieties over $\CC$. Let $H$ be a finite-dimensional complex
vector space. A {\em pure Hodge structure of weight $k$} on $H$ is
a decomposition
 $$
 H=\bigoplus\limits_{p+q=k} H^{p,q}
 $$
such that $H^{q,p}=\overline{H}^{p,q}$, the bar denoting complex
conjugation in $H$. We denote
  $$
  h^{p,q}(H)=\dim H^{p,q}\ ,
  $$
which is called the Hogde number of type $(p,q)$. A Hodge
structure of weight $k$ on $H$ gives rise to the so-called {\em
Hodge filtration} $F$ on $H$, where
  $$
  F^p= \bigoplus\limits_{s\geq p} H^{s,p-s}\ ,
  $$
which is a descending filtration. Note that $\Gr_F^p H =
F^p/F^{p+1}= H^{p,q}$.

Let $H$ be a finite-dimensional complex vector space. A {\em
(mixed) Hodge structure} over $H$ consists of an ascending weight
filtration $W$ on $H$ and a descending Hodge filtration $F$ on $H$
such that $F$ induces a pure Hodge filtration of weight $k$ on
each $\Gr^W_k H= W_k/W_{k-1}$. Again we define
  $$
  h^{p,q}(H)=\dim \, H^{p,q}\, , \qquad \text{where}\quad H^{p,q}
  =\Gr_F^p\Gr^W_{p+q} H\, .
  $$

Deligne has shown \cite{De} that, for each complex algebraic variety
$Z$,  the cohomology $H^k(Z)$ and the cohomology with compact
support $H_c^k(Z)$ both carry natural Hodge structures. If $Z$ is a
compact smooth projective variety (hence compact K{\"a}hler) then
the Hodge structure $H^k(Z)$ is pure of weight $k$ and coincides
with the classical Hodge structure given by the Hodge decomposition
of harmonic forms into $(p,q)$ types.

\begin{definition} \label{def:Hodge-poly}
For \emph{any} complex algebraic variety $Z$ (not necessarily
smooth, compact or irreducible), we define the Hodge numbers as
  $$
  h^{k,p,q}_c(Z)=h^{p,q}(H^k_c(Z))=\dim \Gr^p_F \Gr_{p+q}^W H^k_c(Z)\, .
  $$
Introduce the Euler characteristic
 $$
   \chi^{p,q}_c(Z) = \sum_k (-1)^k h^{k,p,q}_c(Z)
 $$
The \emph{Hodge polynomial} of $Z$ is defined \cite{DK} as
 $$
 e(Z)=e(Z)(u,v)= \sum_{p,q} (-1)^{p+q}\chi_c^{p,q}(Z) u^p v^q\, .
 $$
\end{definition}

If $Z$ is smooth and projective then the mixed Hodge structure on
$H^k_c(Z)$ is pure of weight $k$, so $\Gr_k^W H_c^k(Z)
=H_c^k(Z)=H^k(Z)$ and the other pieces $\Gr_m^W H_c^k(Z)=0$, $m\neq
k$. So
   $$
   \chi_c^{p,q}(Z)=(-1)^{p+q} h^{p,q}(Z),
   $$
where $h^{p,q}(Z)$ is the usual Hodge number of $Z$. In this case,
  $$
  e(Z)(u,v)= \sum_{p,q} h^{p,q}(Z) u^p v^q\,
  $$
is the (usual) Hodge polynomial of $Z$. Note that in this case,
the Poincar{\'e} polynomial of $Z$ is
  \begin{equation}\label{eqn:Poinca}
  P_Z(t)=\sum_k b^k(Z) t^k= \sum_k \left(
  \sum_{p+q=k} h^{p,q}(Z) \right) t^k= e(Z)(t,t).
  \end{equation}
where $b^k(Z)$ is the $k$-th Betti number of $Z$.

\begin{theorem}[{\cite[Theorem 2.2]{MOV}}] \label{thm:Du}
 Let $Z$ be a complex algebraic variety. Suppose that $Z$ is
 a finite disjoint union $Z=Z_1\cup \cdots \cup Z_n$, where the
 $Z_i$ are algebraic subvarieties. Then
  $$
  e(Z)= \sum_i e (Z_i).
  $$\hfill $\Box$
\end{theorem}

Note that we can assign to \textit{any} complex algebraic variety
$Z$ (not necessarily smooth, compact or irreducible) a polynomial
 $$
 P_Z(t) = e(Z)(t,t)= \sum_m (-1)^m \chi^m_c(Z) \, t^m=\sum_{k,m} (-1)^{k+m}
  \dim \Gr_m^W H_c^k(Z) \, t^m ,
 $$
where
 $$
 \chi^m_c (Z)= \sum_{p+q=m}  \chi^{p,q}_c(Z).
 $$
This is called the {\em virtual Poincar{\'e} polynomial} of $Z$ (see
\cite{FM,Du}). It satisfies an additive property analogous to that
of Theorem \ref{thm:Du} and it recovers the usual Poincar{\'e}
polynomial when $Z$ is a smooth projective variety.

The following Hodge polynomials will be needed later:
\begin{itemize}
 \item Let $Z=\PP^{n}$, then $e(Z)= 1+uv+(uv)^2+\cdots +
 (uv)^{n}=(1-(uv)^{n+1})/(1-uv)$. For future reference,
 we shall denote
   \begin{equation}\label{eqn:Pn}
    e_n := e
    (\PP^{n-1})=e
    (\PP(\CC^n))
   =\frac{1-(uv)^{n}}{1-uv}\ .
   \end{equation}
 \item Let $\Jac^d X$ be the Jacobian of (any) degree $d$ of a
 (smooth, projective) complex curve $X$ of genus $g$. Then
  \begin{equation}\label{eqn:Jac}
  e(\Jac^d X)=(1+u)^g(1+v)^g.
   \end{equation}
  \end{itemize}

\begin{lemma}[{\cite[Lemma 2.3]{MOV}}]\label{lem:vb}
  Suppose that  $\pi:Z\to Y$ is an algebraic fiber bundle with fiber $F$ which is
 locally trivial in the Zariski topology, then $e(Z)=e(F)\,e(Y)$.
 (In particular this is true for $Z=F\times Y$.) \hfill $\Box$
\end{lemma}

\begin{lemma}\label{lem:Gri}
  Suppose that $\pi:Z\to Y$ is a map between quasi-projective varieties which is
  a locally trivial fiber bundle in the usual topology, with fibers projective
  spaces $F=\PP^N$ for some $N>0$. Then $e(Z)=e(F)\,e(Y)$.
\end{lemma}

\begin{proof}
This follows from \cite{Gri,Del}. For completeness we provide a
proof. Let $H$ be a hyperplane section of $Z$. We have a morphism of
Hodge structures:
 \begin{equation}\label{eqn:Leray}
 \begin{array}{ccc}
  L: H^*(\PP^N)\ox H^*_c(Y) &\to& H^*_c(Z) \\
  h^i \ox \a &\mapsto & H^i\cap \pi^*(\a)\, ,
 \end{array}
 \end{equation}
where $h$ is the hyperplane class of the projective space. Note that
$L$ is not multiplicative. Let us see that $L$ is injective. If
$x=\sum H^i \cap \pi^*(\a_i)=0$, let $i_0$ be the maximum $i$ for
which $\a_{i}\neq 0$. Then
 $$
 0 = \pi_*(H^{N-i_0}\cap x)=\a_{i_0}\, .
 $$
So $L$ must be injective. On the other hand, the Leray spectral
sequence of the fibration $\pi$ has $E_2$-term isomorphic to
$H^*(\PP^N)\ox H^*_c(Y)$ and converges to $H^*_c(Z)$. So $\dim
H^*(\PP^N)\ox H^*_c(Y) \geq \dim H^*_c(Z)$ and $L$ must be bijective.
Therefore $L$ is an isomorphism of Hodge structures, and the result
follows.
\end{proof}

\begin{lemma}\label{lem:Gr}
The Hodge polynomial of the Grassmannian $\Gr(k,N)$ is
 $$
  e(\Gr(k,N))=
   \frac{(1- (uv)^{N-k+1}) \cdots (1- (uv)^{N-1})
   (1- (uv)^{N})}{(1- uv)\cdots(1-(uv)^{k-1}) (1-(uv)^{k})}\, .
 $$
\end{lemma}

\begin{proof}
This is well-known, but we provide a proof for completeness.

Let us review first the case of the projective space
$\PP^{N-1}=(\CC^N -\{0\}) /(\CC -\{0\})$. Then $\CC^N -\{0\} \to
\PP^{N-1}$ is a locally trivial fibration, since it is the
restriction of the universal line bundle $U \to \PP^{N-1}$ to the
complement of the zero section. Using either Lemma \ref{lem:vb} or
Lemma \ref{lem:Gri}, we have $e(\CC^N-\{0\})=e(\CC-\{0\})\,
e(\PP^{N-1})$, i.e.\ $(uv)^{N}-1=(uv-1) e(\PP^{N-1})$, from where
(\ref{eqn:Pn}) is recovered. Now in the case of $k>1$, denote
 $$
 F(k,n)= \{ (v_1,\ldots,v_k) \, | \, \text{$v_i$ are linearly
 independent vectors of $\CC^n$}\}\ .
 $$
Then $\Gr(k,N)=F(k,N)/\GL(k,\CC)$ and there is a locally trivial
fibration $F(k,N) \to \Gr(k,N)$ with fiber $\GL(k,\CC)\cong F(k,k)$
(again it is the principal bundle associated to the universal bundle
$U\to \Gr(k,N)$). So by Lemma \ref{lem:vb},
$e(\Gr(k,N))=e(F(k,N))/e(F(k,k))$. Now we use that the map
 $$
 F(k,n) \longrightarrow F(k-1,n),
 $$
given by forgetting the last vector, is a locally trivial fibration,
with fiber $\CC^n-\CC^{k-1}$. Using Lemma \ref{lem:vb} and Theorem
\ref{thm:Du}, we have $e(F(k,n))=e(F(k-1,n))\, e(\CC^n-\CC^{k-1})=
e(F(k-1,n))\, ((uv)^{n}-(uv)^{k-1})$. By recursion this gives
 $$
 e(F(k,n))=((uv)^{n}-(uv)^{k-1}) \cdots ((uv)^{n}-uv) ((uv)^{n}-1)\, .
 $$
So
 \begin{eqnarray*}
  e(\Gr(k,N)) &= & \frac{((uv)^{N}-(uv)^{k-1}) \cdots ((uv)^{N}-uv)((uv)^{N}-1)}{((uv)^{k}-(uv)^{k-1})
   \cdots  ((uv)^{k}-uv)((uv)^{k}-1)} \\ &=&
   \frac{(1- (uv)^{N-k+1}) \cdots (1- (uv)^{N-1}) (1- (uv)^{N})}{(1- uv)\cdots(1-(uv)^{k-1}) (1-(uv)^{k})
   } \, .
 \end{eqnarray*}
\end{proof}

\begin{lemma} \label{lem:Z2}
 Let $M$ be a smooth projective variety. Consider the algebraic variety
 $Z=(M\x M)/\ZZ_2$, where $\ZZ_2$ acts as $(x,y)\mapsto (y,x)$. The Hodge
 polynomial of $Z$ is
   $$
   e(Z)= \frac12 \Big(e(M)(u,v)^2 + e(M)(-u^2,-v^2)\Big) \, .
   $$
\end{lemma}

\begin{proof}
The cohomology of $Z$ is
 $$
 H^*(Z)=H^*(M\x M)^{\ZZ_2}=(H^*(M)\ox H^*(M))^{\ZZ_2}\, .
 $$
This is an equality of Hodge structures. The Hodge structure of
$M$ is of pure type, therefore the Hodge structure of $Z$ is also
of pure type. Moreover,
 $$
 H^{p,q}(Z) = \left( \bigoplus_{p_1+p_2=p \atop
q_1+q_2=q}  H^{p_1,q_1}(M)\ox H^{p_2,q_2}(M) \right)^{\ZZ_2}\, .
 $$
Therefore we have
 $$
 h^{p,q}(Z) =\frac12 \hspace{-4mm} \sum_{{p_1+p_2=p \atop q_1+q_2=q} \atop (p_1,q_1)\neq
 (p_2,q_2)} \hspace{-5mm}
  h^{p_1,q_1}(M)h^{p_2,q_2}(M) \quad + \  \epsilon_{p,q}\, ,
  $$
where
 $$
 \epsilon_{p,q}=
 \left\{\begin{array}{ll} 0, & \text{$p$ or $q$ odd}, \\
 \dim \left(\Sym^2 H^{p_1,q_1}(M)\right) , \qquad & p=2p_1, q=2q_1, p_1+q_1 \text{ even}, \\
 \dim \left(\bigwedge^2 H^{p_1,q_1}(M)\right), & p=2p_1, q=2q_1, p_1+q_1 \text{
 odd}.
 \end{array}\right.
 $$
If $V$ is a vector space of dimension $n$, then $\dim \left(\Sym^2
V\right)=\frac12 (n^2+n)$ and $\dim \left(\bigwedge^2
V\right)=\frac12 (n^2-n)$, so
 $$
 \epsilon_{p,q}=\left\{\begin{array}{ll} 0, & \text{$p$ or $q$ odd}, \\
 \frac12 ( h^{p_1,q_1}(M)^2 + (-1)^{p_1+q_1} h^{p_1,q_1}(M) ),\qquad
  & p=2p_1, \, q=2q_1\, . \end{array}\right.
 $$

This yields
 \begin{eqnarray*}
 e(Z) &=& \sum h^{p,q}(Z)u^pv^q \\ &=& \frac12 \sum
 h^{p_1,q_1}(M)h^{p_2,q_2}(M)u^{p_1+p_2}v^{q_1+q_2} + \frac12
 \sum (-1)^{p_1+q_1} h^{p_1,q_1}(M) u^{2p_1}v^{2q_1} \\ &=& \frac12 \, e(M)\cdot
 e(M) + \frac12 \, e(M)(-u^2,-v^2)\, .
 \end{eqnarray*}
\end{proof}

\section{Moduli spaces of triples}
\label{sec:stable-triples}

\subsection{Holomorphic triples}
\label{subsec:triples-definitions}

Let $X$ be a smooth projective curve of genus $g\geq 2$ over $\CC$.
A \emph{holomorphic triple} $T = (E_{1},E_{2},\phi)$ on $X$ consists
of two holomorphic vector bundles $E_{1}$ and $E_{2}$ over $X$, of
ranks $n_1$ and $n_2$ and degrees $d_1$ and $d_2$, respectively, and
a holomorphic map $\phi \colon E_{2} \to E_{1}$. We refer to
$(n_1,n_2,d_1,d_2)$ as the \emph{type} of $T$, to $(n_1,n_2)$ as the
\emph{rank} of $T$, and to $(d_1,d_2)$ as the \emph{degree} of $T$.

A homomorphism from $T' = (E_1',E_2',\phi')$ to $T =
(E_1,E_2,\phi)$ is a commutative diagram
  \begin{displaymath}
  \begin{CD}
    E_2' @>\phi'>> E_1' \\
    @VVV @VVV  \\
    E_2 @>\phi>> E_1,
  \end{CD}
  \end{displaymath}
where the vertical arrows are holomorphic maps. A triple
$T'=(E_1',E_2',\phi')$ is a subtriple of $T = (E_1,E_2,\phi)$ if
$E_1'\subset E_1$ and $E_2'\subset E_2$ are subbundles,
$\phi(E_2')\subset E_1'$ and $\phi'=\phi|_{E_2'}$. A subtriple
$T'\subset T$ is called \emph{proper} if $T'\neq 0 $ and $T'\neq
T$. The quotient triple $T''=T/T'$ is given by $E_1''=E_1/E_1'$,
$E_2''=E_2/E_2'$ and $\phi'' \colon E_2''\to E_1''$ being the map
induced by $\phi$. We usually denote by $(n_1',n_2',d_1',d_2')$
and $(n_1'',n_2'',d_1'',d_2'')$, the types of the subtriple $T'$
and the quotient triple $T''$.

\begin{definition} \label{def:s-slope}
For any $\s \in \RR$ the \emph{$\s$-slope} of $T$ is defined by
 $$
   \mu_{\s}(T)  =
   \frac{d_1+d_2}{n_1+n_2} + \s \frac{n_{2}}{n_{1}+n_{2}}\ .
 $$
To shorten the notation, we define the \emph{$\mu$-slope} and
\emph{$\lambda$-slope} of the triple $T$ as $\mu=\mu(E_{1} \oplus
E_{2})= \frac{d_1+d_2}{n_1+n_2}$ and
$\lambda=\frac{n_{2}}{n_{1}+n_{2}}$, so that $\mu_{\s}(T)=\mu+\s
\lambda$.
\end{definition}

\medskip

\begin{definition}\label{def:sigma-stable}
We say that a triple $T = (E_{1},E_{2},\phi)$ is
\emph{$\s$-stable} if
  $$
  \mu_{\s}(T') < \mu_{\s}(T) ,
  $$
for any proper subtriple $T' = (E_{1}',E_{2}',\phi')$. We define
\emph{$\s$-semistability} by replacing the above strict inequality
with a weak inequality. A triple is called \emph{$\s$-polystable}
if it is the direct sum of $\s$-stable triples of the same
$\s$-slope. It is \emph{$\s$-unstable} if it is not
$\s$-semistable, and \emph{strictly $\s$-semistable} if it is
$\s$-semistable but not $\s$-stable. A $\s$-destabilizing
subtriple $T'\subset T$ is a proper subtriple satisfying
$\mu_{\s}(T') \geq \mu_{\s}(T)$.
\end{definition}

We denote by
  $$
  \cN_\s = \cN_\s(n_1,n_2,d_1,d_2)
  $$
the moduli space of $\s$-polystable triples $T =
(E_{1},E_{2},\phi)$ of type $(n_1,n_2,d_1,d_2)$, and drop the type
from the notation when it is clear from the context. The open
subset of $\s$-stable triples is denoted by $\cN_\s^s =
\cN_\s^s(n_1,n_2,d_1,d_2)$. This moduli space is constructed in
\cite{BGP} by using dimensional reduction. A direct construction
is given by Schmitt \cite{Sch} using geometric invariant theory.

\bigskip

There are certain necessary conditions in order for
$\s$-semistable triples to exist. Let $\mu_i=\mu(E_i)=d_i/n_i$
stand for the slope of $E_i$, for $i=1,2$. We write
  \begin{align*}
  \s_m = &\mu_1-\mu_2\ ,  \\[5pt]
  \s_M = & \left(1+ \frac{n_1+n_2}{|n_1 - n_2|}\right)(\mu_1 - \mu_2)\ ,
      \qquad \mbox{if $n_1\neq n_2$\ .}
  \end{align*}

\begin{proposition}\cite{BGPG} \label{prop:alpha-range}
The moduli space $\cN_\s(n_1,n_2,d_1,d_2)$ is a complex projective
variety. For $n_1, n_2>0$, let $I$ denote the interval
$I=[\s_m,\s_M]$ if $n_1\neq n_2$, or $I=[\s_m,\infty)$ if $n_1=n_2$.
A necessary condition for $\cN_\s(n_1,n_2,d_1,d_2)$ to be non-empty
is that $\s\in I$. \hfill $\Box$
\end{proposition}

\bigskip

\subsection{Critical values}\label{subsec:critical-values}

To study the dependence of the moduli spaces $\cN_\s$ on the
parameter, we need to introduce the concept of critical value
\cite{BGP,MOV}.

\begin{definition}\label{def:critical}
The values of $\s_c\in I$ for which there exist $0 \le n'_1 \leq
n_1$, $0 \le n'_2 \leq n_2$, $d'_1$ and $d'_2$, with $n_1'n_2\neq
n_1n_2'$, such that
 \begin{equation}\label{eqn:sigmac}
 \s_c=\frac{(n_1+n_2)(d_1'+d_2')-(n_1'+n_2')(d_1+d_2)}{n_1'n_2-n_1n_2'},
 \end{equation}
are called \emph{critical values}.
\end{definition}

Given a triple $T=(E_1,E_2,\phi)$, the condition of
$\s$-(semi)stability for $T$ can only change when $\s$ crosses a
critical value. If $\s=\s_c$ as in (\ref{eqn:sigmac}) and if $T$
has a subtriple $T'\subset T$ of type $(n_1',n_2',d_1',d_2')$,
then $\mu_{\s_c}(T')=\mu_{\s_c}(T)$ and
 \begin{enumerate}\itemsep=5pt
 \item if $\lambda'>\lambda$ (where $\lambda'$ is the $\lambda$-slope of
 $T'$), then $T$ is not $\s$-stable for
 $\s>\s_c$,
\item if $\lambda'<\lambda$, then $T$ is not $\s$-stable for
 $\s<\s_c$.
 \end{enumerate}
Note that $n_1'n_2\neq n_1n_2'$ is equivalent to $\lambda'\neq
\lambda$.

\bigskip

Of course, it may happen that there is no triple $T$ as above and hence that the
moduli spaces $\cN_\s$ and $\cN_\s^s$ do not change when crossing $\s_c$ (see Remark
\ref{rem:virtual-not-all}).

\medskip

\begin{proposition}[{\cite[Proposition 2.6]{BGPG}}] \label{prop:triples-critical-range}
Fix $(n_1,n_2,d_1,d_2)$. Then
 \begin{enumerate}
 \item[(1)] The critical values are a finite number of values $\s_c \in I$.
 \item[(2)] The stability and semistability criteria  for two values of $\s$
  lying between two consecutive critical values are equivalent; thus
  the corresponding moduli spaces are isomorphic.
 \item[(3)] If $\s$ is not a critical value and $\GCD(n_1,n_2,d_1+d_2) = 1$,
  then $\s$-semistability is equivalent to $\s$-stability, i.e.,\
  $\cN_\s=\cN_\s^s$.
 \end{enumerate} \hfill $\Box$
\end{proposition}

Note that if $\GCD(n_1,n_2,d_1+d_2)\neq 1$ then it may happen that
there exists triples $T$ which are strictly $\s$-semistable for
non-critical values of $\s$.

\subsection{Extensions and deformations of triples}
\label{subsec:extensions-of-triples}

The homological algebra of triples is controlled by the
hypercohomology of a certain complex of sheaves which appears when
studying infinitesimal deformations \cite[Section 3]{BGPG}. Let
$T'=(E'_1,E'_2,\phi')$ and $T''=(E''_1,E''_2,\phi'')$ be two
triples of types $(n_{1}',n_{2}',d_{1}',d_{2}')$ and
$(n_{1}'',n_{2}'',d_{1}'',d_{2}'')$, respectively. Let
$\Hom(T'',T')$ denote the linear space of homomorphisms from $T''$
to $T'$, and let $\Ext^1(T'',T')$  denote the linear space of
equivalence classes of extensions of the form
 $$
  0 \lto T' \lto T \lto T'' \lto 0,
 $$
where by this we mean a commutative  diagram
  $$
  \begin{CD}
  0@>>>E_1'@>>>E_1@>>> E_1''@>>>0\\
  @.@A\phi' AA@A \phi AA@A \phi'' AA\\
  0@>>>E'_2@>>>E_2@>>>E_2''@>>>0.
  \end{CD}
  $$
To analyze $\Ext^1(T'',T')$ one considers the complex of sheaves
 \begin{equation} \label{eqn:extension-complex}
    C^{\bullet}(T'',T') \colon ({E_{1}''}^{*} \otimes E_{1}') \oplus
  ({E_{2}''}^{*} \otimes E_{2}')
  \overset{c}{\lto}
  {E_{2}''}^{*} \otimes E_{1}',
 \end{equation}
where the map $c$ is defined by
 $$
 c(\psi_{1},\psi_{2}) = \phi'\psi_{2} - \psi_{1}\phi''.
 $$

\begin{proposition}[{\cite[Proposition 3.1]{BGPG}}]
  \label{prop:hyper-equals-hom}
  There are natural isomorphisms
  \begin{align*}
    \Hom(T'',T') &\cong \HH^{0}(C^{\bullet}(T'',T')), \\
    \Ext^{1}(T'',T') &\cong \HH^{1}(C^{\bullet}(T'',T')),
  \end{align*}
and a long exact sequence associated to the complex
$C^{\bullet}(T'',T')$:
 $$
 \begin{array}{c@{\,}c@{\,}c@{\,}l@{\,}c@{\,}c@{\,}c}
  0 &\lto \mathbb{H}^0(C^{\bullet}(T'',T')) &
  \lto & H^0(({E_{1}''}^{*} \otimes E_{1}') \oplus ({E_{2}''}^{*} \otimes
  E_{2}'))
  & \lto &  H^0({E_{2}''}^{*} \otimes E_{1}') \\[3pt]
    &  \lto \mathbb{H}^1(C^{\bullet}(T'',T')) &
  \lto &  H^1(({E_{1}''}^{*} \otimes E_{1}') \oplus ({E_{2}''}^{*} \otimes
  E_{2}'))
 &  \lto & H^1({E_{2}''}^{*} \otimes E_{1}') \\[3pt]
 &   \lto \mathbb{H}^2(C^{\bullet}(T'',T')) & \lto & 0. & &
 \end{array}
 $$ \hfill $\Box$
\end{proposition}

We introduce the following notation:
\begin{align*}
  h^{i}(T'',T') &= \dim\HH^{i}(C^{\bullet}(T'',T')), \\ 
  \chi(T'',T') &= h^0(T'',T') - h^1(T'',T') + h^2(T'',T'). 
\end{align*}

\begin{proposition}[{\cite[Proposition 3.2]{BGPG}}]
  \label{prop:chi(T'',T')}
  For any holomorphic triples $T'$ and $T''$ we have
  \begin{align*}
    \chi(T'',T') &= \chi({E_{1}''}^{*} \otimes E_{1}')
    + \chi({E_{2}''}^{*} \otimes E_{2}')
    - \chi({E_{2}''}^{*} \otimes E_{1}')  \\[5pt]
    &= (1-g)(n''_1 n'_1 + n''_2 n'_2 - n''_2 n'_1) + n''_1 d'_1 - n'_1 d''_1
    + n''_2 d'_2 - n'_2 d''_2
    - n''_2 d'_1 + n'_1 d''_2,
  \end{align*}
where $\chi(E)=\dim H^0(E) - \dim H^1(E)$ is the Euler
characteristic of $E$. \hfill $\Box$
\end{proposition}

\bigskip
Since the  space of infinitesimal deformations of $T$ is
isomorphic to $\HH^{1}(C^{\bullet}(T,T))$, the previous results
also apply to studying deformations of a holomorphic triple $T$.

\begin{theorem}[{\cite[Theorem 3.8]{BGPG}}]\label{thm:smoothdim}
Let $T=(E_1,E_2,\phi)$ be an $\s$-stable triple of type
$(n_1,n_2,d_1,d_2)$.
 \begin{enumerate}
 \item[(1)] The Zariski tangent space at the point defined by $T$
 in the moduli space of stable triples  is isomorphic to
 $\HH^{1}(C^{\bullet}(T,T))$.
 \item[(2)] If\/ $\HH^{2}(C^{\bullet}(T,T))= 0$, then the moduli space of
 $\s$-stable triples is smooth in  a neighbourhood of the point
 defined by $T$.
 \item[(3)] At a smooth point $T\in \cN^s_\s(n_1,n_2,d_1,d_2)$ the
 dimension of the moduli space of $\s$-stable triples is
 \begin{align*}
  \dim \cN^s_\s(n_1,n_2,d_1,d_2)
  &= h^{1}(T,T) = 1 - \chi(T,T) \\ 
  &= (g-1)(n_1^2 + n_2^2 - n_1 n_2) - n_1 d_2 + n_2 d_1 + 1.
 \end{align*}
 \item[(4)] Let $T=(E_1,E_2,\phi)$ be a $\s$-stable triple. If $T$
 is injective or surjective (meaning that $\phi:E_2\to E_1$ is
 injective or surjective) then the moduli space is smooth at $T$.
 \end{enumerate} \hfill $\Box$
\end{theorem}

\subsection{Crossing critical values}
\label{subsec:crossing-critical-values}

Fix the type $(n_1,n_2,d_1,d_2)$ for the moduli spaces of
holomorphic triples. We want to describe the differences between
two spaces $\cN^s_{\s_1}$ and $\cN^s_{\s_2}$ when $\s_1$ and
$\s_2$ are separated by a critical value. Let $\s_c\in I$ be a
critical value and set
 $$
 \scp = \s_c + \epsilon,\quad \scm = \s_c -
 \epsilon,
 $$
where $\epsilon > 0$ is small enough so that $\s_c$ is the only
critical value in the interval $(\scm,\scp)$.

\begin{definition}\label{def:flip-loci}
We define the \textit{flip loci} as
 \begin{align*}
 \cS_{\scp} &= \{ T\in\cN_{\scp} \ |
 \ \text{$T$ is $\scm$-unstable}\} \subset\cN_{\scp} \ ,\\
 \cS_{\scm} &= \{ T\in\cN_{\scm} \ |
 \ \text{$T$ is $\scp$-unstable}\}
 \subset\cN_{\scm} \ .
 \end{align*}
and $\cS_{\scpm}^s=\cS_\scpm \cap \cN_\scpm^s$ for the stable part
of the flip loci.
\end{definition}

Note that for $\s_c=\s_m$, $\cN_{\s_m^-}$ is empty, hence
$\cN_{\smp}= \cS_{\smp}$.  Analogously, when $n_1\neq n_2$,
$\cN_{\s_M^+}$ is empty and $\cN_{\s_M^-}= \cS_{\s_M^-}$.

\begin{lemma}\label{lem:fliploci}
 Let $\s_c$ be a critical value. Then
 \begin{itemize}\itemsep=5pt
 \item[(1)] $\cN_{\scp}-\cS_{\scp}=\cN_{\scm}-\cS_{\scm}$.
 \item[(2)] $\cN^s_{\scp}-\cS_{\scp}^s=
   \cN^s_{\scm}-\cS_{\scm}^s=\cN^s_{\s_c}$.
 \end{itemize}
\end{lemma}

\begin{proof}
Item (1) is an easy consequence of the definition of flip loci.
Item (2) is the content of \cite[Lemma 5.3]{BGPG}.
\end{proof}

Let us describe the flip loci $\cS_{\s_c^{\pm}}$. Let $\s_c$ be a
critical value, and let $(n_1',n_2',d_1',d_2')$ such that
$\lambda' \neq \lambda$ and (\ref{eqn:sigmac}) holds. Put
$(n_1'',n_2'',d_1'',d_2'')=(n_1-n_1',n_2-n_2',d_1-d_1',d_2-d_2')$.
Denote $\cN_\s'=\cN_\s(n_1',n_2',d_1',d_2')$ and
$\cN_\s''=\cN_\s(n_1'',n_2'',d_1'',d_2'')$.

\begin{lemma}[{\cite[Lemma 4.7]{MOV}}] \label{lem:semistable}
Let $T\in \cS_\scp$ (resp.\ $T\in \cS_\scm$). Then $T$ sits in a
non-split exact sequence
 \begin{equation}\label{eqn:extension-triples}
 0\to T'\to T\to T''\to 0,
 \end{equation}
where $\mu_{\s_c}(T')=\mu_{\s_c}(T)=\mu_{\s_c}(T'')$,
$\lambda'<\lambda$ (resp.\ $\lambda'
 >\lambda$) and $T'$ and $T''$ are both
$\s_c$-semistable.

Conversely, if $T'\in \cN_{\s_c}'$ and $T''\in \cN_{\s_c}''$ are
both $\s_c$-stable, and $\lambda'<\lambda$ (resp.\ $\lambda'
>\lambda$). Then for any non-trivial extension {\rm (\ref{eqn:extension-triples})}, $T$
lies in $\cS_\scp^s$ (resp.\ in $\cS_\scm^s$). Moreover, such $T$
can be written uniquely as an extension {\rm
(\ref{eqn:extension-triples})} with
$\mu_{\s_c}(T')=\mu_{\s_c}(T)$.

In particular, suppose $\s_c$ is not a critical value for the
moduli spaces of triples of types $(n_1',n_2',d_1',d_2')$ and
 $(n_1'',n_2'',d_1'',d_2'')$, $\GCD(n_1',n_2',d_1'+d_2')=1$ and
 $\GCD(n_1'',n_2'',d_1''+d_2'')=1$. Then if $\lambda'<\lambda$ (resp.\
 $\lambda'>\lambda$), there is a bijective correspondence
 between non-trivial extensions {\rm (\ref{eqn:extension-triples})},
with
 $T'\in \cN_{\s_c}'$ and $T''\in \cN_{\s_c}''$ and triples
 $T\in \cS_\scp$ (resp.\ $\Smenos$). \hfill $\Box$
\end{lemma}

\medskip

\begin{theorem} \label{thm:Smas}
 Let $\s_c$ be a critical value with $\lambda'<\lambda$
 (resp.\ $\lambda'>\lambda$). Assume
 \begin{itemize}
 \item[(i)] $\s_c$ is not a critical value for the moduli spaces
 of triples of types $(n_1',n_2',d_1',d_2')$ and
 $(n_1'',n_2'',d_1'',d_2'')$, $\GCD(n_1',n_2',d_1'+d_2')=1$ and
 $\GCD(n_1'',n_2'',d_1''+d_2'')=1$.
 \item[(ii)] $\HH^0(C^\bullet(T'',T'))=\HH^2(C^\bullet(T'',T'))
 =0$, for every $(T',T'')\in \cN_{\s_c}'\times  \cN_{\s_c}''$.
 \end{itemize}
 Then $\Smas$ (resp.\ $\Smenos$) is the projectivization of a
 bundle of
 rank $-\chi(T'',T')$ over $\cN_{\s_c}' \times \cN_{\s_c}''$.
\end{theorem}

\begin{proof}
This is the content of \cite[Theorem 4.8]{MOV}. Note that by
\cite{Sch}, the moduli spaces $\cN_{\s_c}'$ and $\cN_{\s_c}''$ are
fine moduli spaces (since $\GCD(n_1',n_2',d_1'+d_2')=1$ and
$\GCD(n_1'',n_2'',d_1''+d_2'')=1$), so the hypothesis (iii) in
\cite[Theorem 4.8]{MOV} is satisfied.
\end{proof}

The construction of the flip loci can be used for the critical
value $\s_c=\s_m$, which allows to describe the moduli space
$\cN_{\smp}$. We refer to the value of $\s$ given by
$\s=\smp=\s_m+\epsilon$ as \textit{small}.

Let $M(n,d)$ denote the moduli space of polystable vector bundles of
rank $n$ and degree $d$ over $X$. This moduli space is projective.
We also denote by $M^s(n,d)$ the open subset of stable bundles,
which is smooth of dimension $n^2(g-1)+1$. If $\GCD(n,d)=1$, then
$M(n,d)=M^s(n,d)$.

\begin{proposition}[{\cite[Proposition 4.10]{MOV}}] \label{prop:moduli-small}
There is a map
 $$
 \pi:\cN_{\smp}=\cN_{\smp}(n_1,n_2,d_1,d_2) \to M(n_1,d_1) \times M(n_2,d_2)
 $$
which sends $T=(E_1,E_2,\phi)$ to $(E_1,E_2)$.
 \begin{enumerate}
 \item[(i)] If $\GCD(n_1,d_1)=1$, $\GCD(n_2,d_2)=1$ and $\mu_1-\mu_2>2g-2$,
then $\cN_{\smp}^s=\cN_{\smp}$ is a projective bundle over
$M(n_1,d_1) \times M(n_2,d_2)$, whose fibers are projective spaces
of dimension $n_2d_1-n_1d_2- n_1n_2(g-1)-1$.
 \item[(ii)] In general, if $\mu_1-\mu_2>2g-2$, then the open subset
 $$
 \pi^{-1}(M^s(n_1,d_1)\times  M^s(n_2,d_2)) \subset \cN_{\smp}
 $$
is a projective bundle over $M^s(n_1,d_1) \times M^s(n_2,d_2)$,
whose fibers are projective spaces of dimension $n_2d_1-n_1d_2-
n_1n_2(g-1)-1$.
 \end{enumerate} \hfill $\Box$
\end{proposition}

\section{Hodge polynomials of the moduli spaces of
triples of ranks $(2,1)$ and $(1,2)$} \label{sec:Poincare(2,1)}

\subsection{Moduli of triples of rank $(2,1)$}
\label{subsec:triples(2,1)}

In this section we recall the main results of \cite{MOV}. Let
$\moduli=\moduli(2,1,d_1,d_2)$ denote the moduli space of
$\s$-polystable triples $T=(E_1,E_2,\phi)$ where $E_1$ is a vector
bundle of degree $d_1$ and rank $2$ and $E_2$ is a line bundle of
degree $d_2$. By Proposition \ref{prop:alpha-range}, $\s$ is in
the interval
  $$
    I=[\s_m,\s_M]=
    [\mu_1-\mu_2\,,\,4(\mu_1-\mu_2)]=[d_1/2-d_2,2d_1-4d_2], \qquad\mbox{ where }
    \mu_1-\mu_2\ge0\, .
  $$
Otherwise $\cN_\s$ is empty.

\begin{theorem}[{\cite[Theorem 5.1]{MOV}}]\label{thm:moduli(2,1)}
For $\s\in I$, $\moduli$ is a projective variety. It is smooth and
of (complex) dimension $3g-2 + d_1 - 2 d_2$ at the stable points
$\moduli^s$. Moreover, for non-critical values of $\s$,
$\moduli=\moduli^s$ (hence it is smooth and projective). \hfill
$\Box$
\end{theorem}

The critical values corresponding to $n_1=2$, $n_2=1$ are given by
Definition \ref{def:critical}\ :

\begin{enumerate}
 \item[(1)] $n'_1=1$, $n'_2=0$. The corresponding $\s_c$-destabilizing
  subtriple is of the form $0\to E_1'$, where $E_1'=M$ is a line bundle of
  degree $\deg(M)=d_M$. The critical value is
  $$
   \s_c= 3d_M-d_1-d_2\, .
  $$
 \item[(2)] $n'_1=1$, $n'_2=1$. The corresponding
 $\s_c$-destabilizing subtriple $T'$ is of the form $E_2\to
 E_1'$, where $E_1'$ is a line bundle. Let $T''=T/T'$ be the
 quotient bundle, which is of the form $0\to E_1''$, where
 $E_1''=M$ is a line bundle, and let $d_M=\deg(M)$ be its degree.
 Then $d_2'=d_2$, $d_1'=d_1-d_M$ and
  $$
   \s_c=-\big( 3(d_1-d_M+d_2)-2(d_1+d_2)\big)= 3d_M-d_1-d_2\, .
  $$

 \item[(3)] $n'_1=2$, $n'_2=0$. In this case, the only possible subtriple
  is $0\to E_1$. This produces the critical value
  $$
   \s_c=\frac{d_1-2d_2}{2}= \mu_1-\mu_2=\s_m\, ,
  $$
  i.e.,\ the minimum of the interval $I$ for
  $\s$.
 \item[(4)] $n'_1=0$, $n'_2=1$. The subtriple $T'$ must be of the form $E_2\to
  0$. This forces $\phi=0$ in $T=(E_1,E_2,\phi)$. So $T$ is decomposable,
  of the form $T'\oplus T''=(0,E_2,0)\oplus (E_1,0,0)$, and $T$ is
  $\s$-unstable for any $\s\neq \s_c$, where
  $$
   \s_c=\frac{2d_2-d_1}{-2}= \mu_1-\mu_2=\s_m\, .
  $$
\end{enumerate}

\begin{lemma}[{\cite[Lemma 5.3]{MOV}}] \label{lem:dM}
Let $\s_c=3d_M-d_1-d_2$ be a critical value. Then
 \begin{equation} \label{eqn:dM-bounds}
  \mu_1\leq d_M\leq d_1-d_2 \, ,
 \end{equation}
and $\s_c=\s_m \Leftrightarrow d_M=\mu_1$. \hfill $\Box$
\end{lemma}

The Hodge polynomials of the moduli spaces $\cN_\s$ for non-critical
values of $\s$ are given in \cite[Theorem 6.2]{MOV}. As this moduli
space is projective and smooth, we may recover the Poincar{\'e}
polynomial from the Hodge polynomial via the formula
(\ref{eqn:Poinca}).

\begin{theorem}[{\cite[Theorem 6.2]{MOV}}] \label{thm:polinomiono(2,1)no-critico}
Suppose that $\s>\s_m$ is not a critical value. Set
$d_0=\Big[\frac13(\s+d_1+d_2)\Big]+1$. Then the Hodge polynomial of
$\moduli=\moduli (2,1,d_1,d_2)$ is
 $$
  e(\moduli)= \coeff_{x^0}
  \left[\frac{(1+u)^{2g}(1+v)^{2g}(1+ux)^{g}(1+vx)^{g}}{(1-uv)(1-x)(1-uvx)x^{d_1-d_2-d_0}}
    \Bigg(\frac{(uv)^{d_1-d_2-d_0}}{1-(uv)^{-1}x}-
    \frac{(uv)^{-d_1+g-1+2d_0}}{1-(uv)^2x}\Bigg)\right] .
 $$ \hfill $\Box$
\end{theorem}

\subsection{Moduli space of triples of rank $(1,2)$} \label{sec:Poincare(1,2)}

Triples of rank $(1,2)$ are of the form $\phi: E_2\to E_1$, where
$E_2$ is a rank $2$ bundle and $E_1$ is a line bundle. By
Proposition \ref{prop:alpha-range}, $\s$ is in the interval
  $$
    I=[\s_m,\s_M]=
    [\mu_1-\mu_2\,,\,4(\mu_1-\mu_2)]=[d_1-d_2/2,4d_1-2d_2], \qquad\mbox{ where }
    \mu_1-\mu_2\ge0\, .
  $$

\begin{theorem}\label{thm:moduli(1,2)}
For $\s\in I$, $\moduli$ is a projective variety. It is smooth and
of (complex) dimension $3g -2 + 2 d_1 - d_2$ at the stable points
$\moduli^s$. Moreover, for non-critical $\s$, $\moduli=\moduli^s$
(hence it is smooth and projective).
\end{theorem}

\begin{proof}
Given a triple  $T=(E_1,E_2,\phi)$ one has the dual triple
$T^*=(E_2^*,E_1^*,\phi^*)$, where $E_i^*$ is the dual of $E_i$ and
$\phi^*$ is the transpose of $\phi$. The map $T\mapsto T^*$ defines
an isomorphism
 $$
 \cN_\s(1,2,d_1,d_2) \cong \cN_\s(2,1,-d_2,-d_1)\, .
 $$
The result now follows from Theorem \ref{thm:moduli(2,1)}.
\end{proof}

Also from Lemma \ref{lem:dM}, we get

\begin{lemma}[{\cite[Lemma 7.2]{MOV}}] \label{lem:dM-2}
The critical values for $\moduli(1,2,d_1,d_2)$ are the numbers
$\s_c=3d_M+d_1+d_2$, where $-\mu_2\leq d_M \leq d_1-d_2$. Also
$\s_c=\s_m \Leftrightarrow d_M=-\mu_2$. \hfill $\Box$
\end{lemma}

\begin{theorem}[{\cite[Theorem 7.3]{MOV}}] \label{thm:polinomiono(1,2)no-critico}
Consider $\cN_\s=\cN_\s(1,2,d_1,d_2)$. Let $\s>\s_m$ be a
non-critical value. Set $d_0=\Big[\frac13(\s-d_1-d_2)\Big]+1$. Then
the Hodge polynomial of $\moduli$ is
 $$
  e(\moduli)= \coeff_{x^0}
 \left[\frac{(1+u)^{2g}(1+v)^{2g}(1+ux)^{g}(1+vx)^{g}}{(1-uv)(1-x)(1-uvx)x^{d_1-d_2-d_0}}
    \Bigg(\frac{(uv)^{d_1-d_2-d_0}}{1-(uv)^{-1}x}-
    \frac{(uv)^{d_2+g-1+2d_0}}{1-(uv)^2x}\Bigg)\right] .
 $$
\end{theorem}

\begin{proof}
We use that $e(\moduli(1,2,d_1,d_2))=e(\moduli(2,1,-d_2,-d_1))$
and the formula in Theorem \ref{thm:polinomiono(2,1)no-critico},
where $d_1$ and $d_2$ are substituted by $-d_2$, $-d_1$ and
 $$
 d_0= \left[ \frac13 (\s -d_2-d_1)\right] +1 \, .
 $$
\end{proof}

\section{Hodge polynomial of the moduli space of
rank $2$ even degree stable bundles} \label{sec:poly-bundles-rk-two}

Let $M(2,d)$ denote the moduli space of polystable vector bundles of
rank $2$ and degree $d$ over $X$. As $M(2,d)\cong M(2,d+2k)$, for
any integer $k$, there are two moduli spaces, depending on whether
the degree is even or odd. We are going to apply the results of the
Section \ref{sec:Poincare(2,1)} to compute the Hodge polynomials of
these moduli spaces.

We first recall the Hodge polynomial of the moduli space of rank
$2$ odd degree stable bundles from \cite{BaR,EK,MOV}

\begin{theorem}[{\cite[Proposition 8.1]{MOV}}]\label{thm:rank2odd}
The Hodge polynomial of $M(2,d)$ with odd degree $d$, is
 $$
    e(M(2,d)) = \frac{(1+u)^{g}(1+v)^g(1+u^2v)^{g}(1+uv^2)^g
    -(uv)^{g}(1+u)^{2g}(1+v)^{2g}}
    {(1-uv)(1-(uv)^2)}\,.
 $$ \hfill $\Box$
\end{theorem}

Now we compute the Hodge polynomial of the moduli space of rank
$2$ even degree stable bundles. Note that this moduli space is
smooth but non-compact. It is irreducible and of dimension $4g-3$.

\begin{theorem}\label{thm:rank2even}
The Hodge polynomial of $M^s(2,d)$ with even degree $d$, is
 \begin{align*}
   e(M^s(2,d))=\, &
  \frac{1}{2(1-uv)(1-(uv)^2)} \bigg( 2(1+u)^{g}(1+v)^g(1+u^2v)^{g}(1+uv^2)^g   \\
  &-
  (1+u)^{2g}(1+v)^{2g} (1+ 2 u^{g+1} v^{g+1} -u^2v^2) -
  (1-u^2)^g(1-v^2)^g(1-uv)^2  \bigg) \, .
 \end{align*}
\end{theorem}

\begin{proof}
We compute this by relating $M^s(2,d)$ with the moduli space
$\cN_{\smp}=\cN_{\smp}(2,1,d,d_2)$ of triples of rank $(2,1)$ for
small $\sigma$. Choose $(n_1,d_1)=(2,d)$ and $(n_2,d_2)=(1,d_2)$.
If $d_2$ is very negative so that $\mu_1-\mu_2=d/2-d_2>2g-2$ then
Proposition \ref{prop:moduli-small} (ii) applies. We shall choose
the maximum possible value of $d_2$ for this condition to hold,
i.e.\ $d-2d_2=4g-2$.

There is a decomposition $\modulimmas= X_0 \sqcup X_1 \sqcup X_2
\sqcup X_3 \sqcup X_4$ into locally closed algebraic subsets,
defined by the following strata:

\begin{itemize}
\item[(1)] The open subset $X_0\subset \cN_{\smp}$ consists of
those triples of the form $\phi: L\to E$, where $E$ is a stable
rank $2$ bundle of degree $d$, $L$ is a line bundle of degree
$d_2$, and $\phi$ is a non-zero map (defined up to multiplication
by non-zero scalars). Actually, by Proposition
\ref{prop:moduli-small} there is a map
  $$
  \pi: \cN_{\smp} \to  M(2,d)\times \Jac^{d_2} X,
  $$
and $X_0=\pi^{-1}(M^s(2,d)\times \Jac^{d_2}X)$. Proposition
\ref{prop:moduli-small} (ii) says that $X_0$ is a projective
bundle over $M^s(2,d)\times \Jac^{d_2}X$ with fibers isomorphic to
$\PP^{d-2d_2-2g+2-1}=\PP^{2g-1}$. By Lemma \ref{lem:Gri},
  $$
  e(X_0)=e(M^s(2,d)) e(\Jac\, X) e_{2g}\, ,
  $$
where $e_{2g}=e(\PP^{2g-1})$ following the notation in
(\ref{eqn:Pn}).

\item[(2)] The subset $X_1$ parametrizes triples $\phi:L\to E$
where $E$ is a strictly semistable bundle of degree $d$ which sits
as a non-trivial extension
  \begin{equation}\label{eqn:111}
  0\to L_1\to E\to L_2\to 0,
  \end{equation}
with $L_1\not\iso L_2$, $L_1,L_2\in \Jac^{d/2} X$ and $L\in
\Jac^{d_2}X$.

Let $Y_1$ be the family which parametrizes such bundles $E$. For
fixed $L_1,L_2$ with $L_1\not\iso L_2$, the extensions
(\ref{eqn:111}) are determined by $\PP \Ext^1(L_2,L_1)$. As
$L_1,L_2$ are non-isomorphic, $\dim \Ext^1(L_2,L_1)=\dim
H^1(L_1\ox L_2^*)=g-1$, so $\PP \Ext^1(L_2,L_1) \cong \PP^{g-2}$.
Therefore $Y_1$ is a fiber bundle over $\Jac^{d/2}X\x \Jac^{d/2}X
 -  \Delta$, where $\Delta$ is the diagonal, with fibers
isomorphic to $\PP^{g-2}$. Thus using Theorem \ref{thm:Du} and Lemma
\ref{lem:Gri},
 \begin{equation}\label{eqn:Y1}
  {}\qquad e(Y_1)= \big(e(\Jac X)^2 - e(\Jac X)\big) e_{g-1} \, .
 \end{equation}

Now we want to describe $X_1$. For each fixed $E\in Y_1$ as in
(\ref{eqn:111}), and $L\in \Jac^{d_2}X$, there is an exact
sequence
 $$
  0\to \Hom(L,L_1) \to \Hom (L,E)\to \Hom(L,L_2) \to 0.
 $$
Here $\Ext^1(L,L_1)=0$ since $\deg(L_1)-\deg(L)=d/2-d_2>2g-2$. So we
may write $\Hom(L,E) \cong \Hom(L,L_1)\oplus \Hom(L,L_2)$,
non-canonically. Let us see when $\phi\in \Hom(L,E)$ gives rise to a
$\smp$-stable triple $T=(E,L,\phi)$. First note that $T$ is
$\s_m$-semistable, since by Section \ref{subsec:triples(2,1)}, the
only possibility for not being $\s_m$-semistable is to have a
subtriple of rank $(0,1)$, i.e., a line subbundle $M\subset E$,
which by Lemma \ref{lem:dM} should have degree $d_M>\mu_1$,
contradicting the semistability of $E$. If $T$ is not $\smp$-stable
then it must have a $\s_m$-destabilizing subtriple $T'$ of rank
$(1,1)$ by Section \ref{subsec:triples(2,1)}. Such subtriple is of
the form $\phi:L\to L'$, with $L'\subset E$. As
$\mu_{\s_m}(T')=\mu_{\s_m}(T) \implies \mu(L')=\mu(E)$, $L'$ is a
destabilizing subbundle of $E$. But the only destabilizing subbundle
of $E$ is $L_1$, so $\phi$ satisfies $\phi(L)\subset L_1$.
Equivalently, $\phi=(\phi_1,0)\in \Hom(L,E)$ gives rise to
$\smp$-unstable triples.

This discussion implies that given $(E,L)\in Y_1 \x \Jac^{d_2} X$,
the morphisms $\phi$ giving rise to $\smp$-stable triples
$(E,L,\phi)$ are those in
 \begin{equation}\label{eqn:X1}
 \Hom(L,E) - \Hom (L,L_1).
 \end{equation}
By Riemann-Roch, $\dim\Hom(L,E)=d-2d_2-2g+2=2g$ and $\dim \Hom
(L,L_1)=d/2-d_2-g+1=g$. So the space (\ref{eqn:X1}) is isomorphic
to $\CC^{2g}-\CC^{g}$.

The isomorphism class of the triple $T=(E,L,\phi)$ is determined
up to multiplication by non-zero scalar $(E,L,\phi)\mapsto (E,L,
\lambda\phi)$, since $\Aut(T)=\CC^*$. This follows from the fact
that $\Aut(E)=\CC^*$ (since $E$ is a non-trivial extension
(\ref{eqn:111})) and $\Aut(L)=\CC^*$. Taking into account the
$\CC^*$-action by automorphisms, the fibers of the map $\pi:
X_1\to Y_1 \x \Jac^{d_2} X$ are isomorphic to the projectivization
of (\ref{eqn:X1}), i.e.\ $\PP^{2g-1}-\PP^{g-1}$. Hence
  $$
  {}\qquad e(X_1)= e(\Jac X) e(Y_1) (e_{2g} -e_{g})
   = e(\Jac X)^2 (e (\Jac X)-1) e_{g-1} (e_{2g}-e_{g}) 
   \, .
  $$
(For this, write $X_1=X_1'-X_1''$, where $X_1'$ is a
$\PP^{2g-1}$-bundle over $Y_1$ and $X_1''$ is a $\PP^{g-1}$-bundle
over $Y_1$. By Theorem \ref{thm:Du}, $e(X_1)=e(X_1')-e(X_1'')$. Now
use Lemma \ref{lem:Gri} to compute $e(X_1')$ and $e(X_1'')$.)

\item[(3)] The subset $X_2$ parametrizes triples $\phi:L\to E$
where $E$ is a strictly semistable bundle of degree $d$ which sits
as a non-trivial extension
  $$
  0\to L_1\to E\to L_1\to 0
  $$
with $L_1 \in \Jac^{d/2} X$ and $L\in \Jac^{d_2}X$.

The family $Y_2$ parametrizing such bundles $E$ is a fiber bundle
over $\Jac^{d/2}X$ with fibers $\PP \Ext^1(L_1,L_1) =\PP
H^1(\cO)=\PP^{g-1}$ (actually, this fiber bundle is trivial, so
$Y_2=\Jac^{d/2} X\x \PP^{g-1}$). Thus by Lemma \ref{lem:vb},
   \begin{equation}\label{eqn:Y2}
  {}\qquad e(Y_2)= e(\Jac X) e_{g} \, .
  \end{equation}

For each $L_1\in \Jac^{d/2} X$, there is an exact sequence
 $$
  0\to \Hom(L,L_1) \to \Hom (L,E)\to \Hom(L,L_1) \to 0.
 $$
So we may write $\Hom(L,E)\cong \Hom(L,L_1)\oplus \Hom(L,L_1)$,
non-canonically. In order to describe $X_2$, let us see when a
triple $T=(E,L,\phi)$, with $E\in Y_2$, is $\smp$-stable. As
before, the morphisms $\phi$ giving rise to $\smp$-stable triples
$(E,L,\phi)$ are those in
 \begin{equation}\label{eqn:X2}
 \Hom(L,E) - \Hom (L,L_1) = \Hom(L,L_1) \x (\Hom(L,L_1)- \{0\}) \,
 .
 \end{equation}

For a bundle $E$ in $Y_2$, the automorphism group of $E$ is
$\CC\times \CC^*$, where $\CC\times \CC^*$ acts on $\Hom(L,E)$ by
 $$
 (a, \lambda)\cdot (\phi_1,\phi_2) =(\lambda \phi_1 + a \phi_2,
 \lambda \phi_2).
 $$
Thus for any $(E,L)\in Y_2 \x \Jac^{d_2} X$, the morphisms $\phi$
giving rise to $\smp$-stable triples $(E,L,\phi)$ are parametrized
by
 \begin{equation}\label{eqn:X22}
 (\Hom(L,L_1) \x (\Hom(L,L_1)- \{0\}))/\CC\times \CC^* \, .
 \end{equation}
This is a fiber bundle over $\PP\Hom(L,L_1)=(\Hom(L,L_1)-
\{0\})/\CC^*$ with fibers isomorphic to $\Hom(L,L_1)/\CC\phi_2$
for every $[\phi_2]\in \PP\Hom(L,L_1)$. As
$\dim\Hom(L,E)=d-2d_2-2g+2=2g$ and $\dim \Hom
(L,L_1)=d/2-d_2-g+1=g$, the space (\ref{eqn:X22}) is a
$\CC^{g-1}$-bundle over $\PP^{g-1}$.

Therefore $X_2\to Y_2 \x \Jac^{d_2} X$ is $\CC^{g-1}$-bundle over
a $\PP^{g-1}$-bundle over $Y_2 \x \Jac^{d_2} X$. So
  $$
  {}\qquad e(X_2)= e(\Jac X)e(Y_2) e_g (e_g-e_{g-1}) =
  e(\Jac X)^2 e_g^2 (e_g-e_{g-1}) 
 \, .
  $$
(To apply Lemma \ref{lem:Gri}, we write $X_2\to P$, where $P$ is
the $\PP^{g-1}$-bundle over $Y_2 \x \Jac^{d_2} X$. Then
$X_2=X_2'-X_2''$, where $X_2'$ is a $\PP^{g-1}$-bundle over $P$
and $X_2''$ is a $\PP^{g-2}$-bundle over $P$.)

\item[(4)]  The subset $X_3$ parametrizes triples $\phi:L\to E$
where $E$ is a decomposable bundle of the form $E=L_1\oplus L_2$,
$L_1\not\cong L_2$, $L_1,L_2\in \Jac^{d/2}X$ and $L\in
\Jac^{d_2}X$. The space parametrizing such bundles $E$ is
 \begin{equation}\label{eqn:Y3}
 Y_3=\tilde{Y}_3/\ZZ_2, \qquad \text{where } \
 \tilde{Y}_3 = \Jac^{d/2}X\x \Jac^{d/2}X -\Delta\, ,
 \end{equation}
with $\ZZ_2$ acting by permuting the two factors.

As before, the condition for $\phi\in \Hom(L,E)$ to give rise to a
$\smp$-unstable triple is that there is a subtriple $\phi:L\to L'$
where $\mu(L')=\mu(E)$. There are only two possible such choices for
$L'$, namely $L_1$ and $L_2$. So given $(E,L)\in Y_3\x \Jac^{d_2}
X$, the morphisms $\phi\in \Hom(L,E)=\Hom(L,L_1)\oplus \Hom(L,L_2)$
giving rise to $\smp$-stable triples $(E,L,\phi)$ are those with
both components non-zero, i.e., lying in
 $$
 (\Hom (L,L_1)-\{0\})\times (\Hom(L,L_2)-\{0\}) .
 $$

The automorphisms of $E$ are $\Aut(E)=\CC^*\times \CC^*$,
therefore the map $\phi\in \Hom(L,E)=\Hom(L,L_1)\oplus
\Hom(L,L_2)$ is determined up to the action of $\CC^*\times \CC^*$
on both factors. So $\phi$ are parametrized by
 $$
 \PP\Hom (L,L_1) \times \PP\Hom(L,L_2).
 $$

Let $\tilde{X}_3\to \tilde{Y}_3 \x \Jac^{d_2}X$ be the fiber
bundle with fiber over $(L_1,L_2,L)$ equal to $\PP\Hom (L,L_1)
\times \PP\Hom(L,L_2)$. Then $X_3=\tilde X_3/\ZZ_2$, where $\ZZ_2$
acts by permuting $(\phi_1,\phi_2)\mapsto (\phi_2,\phi_1)$. This
covers the action of $\ZZ_2$ on $\tilde{Y}_3$. Now $\tilde
X_3=X_3'-X_3''$, where $\pi: X_3'\to \Jac^{d/2}X\x \Jac^{d/2}X\x
\Jac^{d_2}X$ is a $\PP^{g-1}\x \PP^{g-1}$-bundle and
$X_3''=\pi^{-1}(\Delta\x \Jac^{d_2}X)$. If $A\to \Jac^{d/2} X\x
\Jac^{d_2}X$ is the $\PP^{g-1}$-bundle with fiber over $L_1$ equal
to $\PP\Hom (L,L_1)$, then $X_3'=A\x_{\Jac^{d_2}X} A$. We apply
Lemma \ref{lem:Z2} fiberwise: $A\to \Jac^{d_2}X$ is a fibration
whose fiber is $A_L$, which in turn is a fibration over $\Jac^{d/2}X$ with
fibers $\PP\Hom (L,L_1)$. Then $X_3'$ fibers over $\Jac^{d_2}X$
with fibers $(A_L\x A_L)/\ZZ_2$. Now
  $$
   \begin{aligned}
 H^*((A\x_{\Jac^{d_2}X} A)/\ZZ_2) &= H^*(A\x_{\Jac^{d_2}X}
 A)^{\ZZ_2} \\ &= (H^*(A_L \x  A_L)\ox H^*(\Jac^{d_2}X))^{\ZZ_2} \\ &=
 (H^*(A_L \x  A_L))^{\ZZ_2}\ox H^*(\Jac^{d_2}X)\, .
  \end{aligned}
 $$
So
  $$
  \begin{aligned}
  e(X_3'/\ZZ_2) & = e((A_L\x A_L)/\ZZ_2) e(\Jac X) \\
  &= \frac12 \Big( e(\Jac X)^2e_g^2 +  (1-u^2)^g(1-v^2)^g
  \frac{1-(uv)^{2g}}{1-u^2v^2}\Big) e(\Jac X)\, .
  \end{aligned}
  $$
On the other hand, $X_3''$ is a $\PP^{g-1}\x \PP^{g-1}$-bundle
over $\Delta\x \Jac^{d_2}X$, the action of $\ZZ_2$ is trivial on
the base, and acts by permutation on the fibers. So $X_3''/\ZZ_2$
is a bundle over $\Delta\x \Jac^{d_2}X$ with fibers
  $$
  (\PP \Hom (L,L_1)\times \PP\Hom(L,L_1))/\ZZ_2
  =(\PP^{g-1}\x \PP^{g-1})/\ZZ_2.
  $$
This fibration is locally trivial in the Zariski topology, since
it is associated to a locally trivial (in the Zariski topology)
vector bundle over $\Delta\x \Jac^{d_2}X$. Hence by Lemma
\ref{lem:vb} and Lemma \ref{lem:Gri},
 $$
  e(X_3''/\ZZ_2)= e(\Jac X)^2 e(\PP^{g-1}\x \PP^{g-1}/\ZZ_2) =
  \frac12e(\Jac X)^2  \bigg( e_g^2 +  \frac{1-(uv)^{2g}}{1-u^2v^2}\bigg)\, .
  $$
Finally using Theorem \ref{thm:Du},
  $$
  \begin{aligned}
  e(X_3)= & \ e(\tilde{X}_3/\ZZ_2)=  e(X_3'/\ZZ_2)-e(X_3''/\ZZ_2) \\
   = & \ \frac12 \bigg( e(\Jac X)^2e_g^2 +  (1-u^2)^g(1-v^2)^g
   \frac{1-(uv)^{2g}}{1-u^2v^2}\bigg) e(\Jac X) 
   -\frac12e(\Jac X)^2  \bigg( e_g^2 +
   \frac{1-(uv)^{2g}}{1-u^2v^2}\bigg)\, .
  \end{aligned}
  $$

\item[(5)] The subset $X_4$ parametrizes triples $\phi:L\to E$,
where $E$ is a decomposable bundle of the form $E=L_1\oplus L_1$,
$L_1\in \Jac^{d/2}X$ and $L\in \Jac^{d_2}X$. Such bundles $E$ are
parametrized by $Y_4=\Jac^{d/2}X$. The morphism $\phi$ lives in
 \begin{equation}\label{eqn:Y4}
 \Hom(L,E)=\Hom(L,L_1)\oplus \Hom(L,L_1)=\Hom(L,L_1)\otimes \CC^2\, .
 \end{equation}

The condition for a triple $T=(E,L,\phi)$ to be $\smp$-unstable is
that there is a destabilizing subbundle $L' \subset E$. A
destabilizing subbundle of $E$ is necessarily isomorphic to $L_1$
and there exists $(a,b)\neq (0,0)$ such that $L'\cong L_1 \inc E$
is given by $x\mapsto (ax, bx)$. This means that
$\phi=(a\psi,b\psi) \in \Hom(L,L_1)\otimes \CC^2$, for some
$\psi\in \Hom(L,L_1)$. All this discussion implies that the set of
$\phi$ giving rise to $\smp$-stable triples are those of the form
$\phi=(\phi_1,\phi_2)\in\Hom(L,L_1)\otimes \CC^2$, with
$\phi_1,\phi_2$ linearly independent.

The automorphisms of $T=(E,L,\phi)$ are $\Aut(T)\cong
\Aut(E)=GL(2,\CC)$. This acts on (\ref{eqn:Y4}) via the standard
representation of $GL(2,\CC)$ on $\CC^2$. So the morphisms $\phi$
are parametrized by the grassmannian $\Gr(2, \Hom(L,L_1))$. As
$\dim \Hom(L,L_1)=g$, we have that $\Gr(2, \Hom(L,L_1))\cong
\Gr(2,g)$.

Moreover $X_4\to Y_4\x \Jac^{d_2}X$ is a locally trivial fibration
in the Zariski topology since it is associated to the (locally
trivial in the Zariski topology) vector bundle over $Y_4\x
\Jac^{d_2}X$ with fibers $\Hom(L,L_1)$. Using Lemma \ref{lem:Gr},
 $$
  e(X_4)= e (\Jac X)^2 e(\Gr(2, g))
  = e (\Jac X)^2
  \frac{(1-(uv)^{g-1})(1-(uv)^{g})}{(1-(uv)^2)(1-uv)}
  \, .
 $$
\end{itemize}

Putting all together,
   \begin{equation}\label{eqn:M(2,even)}
  \begin{aligned}
    e(\modulimmas)=& \
    e(X_0)+e(X_1)+ e(X_2)+e(X_3)+e(X_4)\\
    =& \ e(M^s(2,d)) e(\Jac X)e_{2g} +  e(\Jac X)^2 (e (\Jac X)-1) e_{g-1} (e_{2g}-e_{g})
      + e(\Jac X)^2 e_g^2 (e_g-e_{g-1}) \\
     & +\frac12 \bigg( e(\Jac X)^2e_g^2 +  (1-u^2)^g(1-v^2)^g
   \frac{1-(uv)^{2g}}{1-u^2v^2}\bigg) e(\Jac X)
  -\frac12e(\Jac X)^2  \bigg( e_g^2 +
  \frac{1-(uv)^{2g}}{1-u^2v^2}\bigg)\\
   &+e (\Jac X)^2 e(\Gr(2, g))\,.
  \end{aligned}
  \end{equation}
To compute the left hand side, we use Theorem
\ref{thm:polinomiono(2,1)no-critico} for
$\s=\smp=\mu_1-\mu_2+\epsilon$, $\epsilon>0$ small. It gives
 $$
 d_0=\big[\mbox{$\frac13$}(\mu_1-\mu_2+\varepsilon+2\mu_1+\mu_2)\big]+1=
 [\mu_1]+1=\frac{d}2+1 \,.
 $$
Substuting into the formula for $e(\moduli)$ with $d_1=d/2$ and
$d-2d_2=4g-2$, the Hodge polynomial of $\modulimmas$ equals
  $$
    e(\modulimmas)=
    \coeff_{x^0}\Bigg[
    \frac{(1+u)^{2g}(1+v)^{2g}(1+ux)^{g}(1+vx)^g}
    {(1-uv)(1-x)(1-uvx)x^{2g-2}}
    \Bigg(   \frac{(uv)^{2g-2}}{1-(uv)^{-1}x} -
    \frac{(uv)^{g+1}}{1-(uv)^2x}\Bigg)\Bigg] \, .
 $$
Using the following equality (see the proof of \cite[Proposition
8.1]{MOV})
 $$
    \coeff_{x^{0}}
    \frac{(1+ux)^{g}(1+vx)^g}{(1-ax)(1-bx)(1-cx)x^{2g-2}}=
    \frac{(a+u)^{g}(a+v)^{g}}{(a-b)(a-c)}+
    \frac{(b+u)^{g}(b+v)^{g}}{(b-a)(b-c)}+
    \frac{(c+u)^{g}(c+v)^{g}}{(c-a)(c-b)}\, ,
 $$
one gets the following expression
 $$
\begin{aligned}
    e(\modulimmas) = & \
    \frac{(1+u)^{2g}(1+v)^{2g}}{(1-uv)^2(1-(uv)^2)}
    \big[(1+u^2v)^{g}(1+v^2u)^g(1-(uv)^{2g})+ \\ & \ +(1+u)^{g} (1+v)^g
    \big((uv)^{3g-1}+(uv)^{2g+1}- (uv)^{2g-1}-(uv)^{g+1}\big)
    \big]\, .
 \end{aligned}
 $$

Finally we substitute this into (\ref{eqn:M(2,even)}) to get the
Hodge polynomial $e(M^s(2,d))$ as in the statement.
\end{proof}

\begin{corollary} \label{cor:rank2even-ss} The Hodge polynomial of the moduli
space of polystable rank $2$ even degree $d$ vector bundles is
 \begin{align*}
   e(M(2,d))=\, &
  \frac{1}{2(1-uv)(1-(uv)^2)} \bigg( 2(1+u)^{g}(1+v)^g(1+u^2v)^{g}(1+uv^2)^g   \\
  &-
  (1+u)^{2g}(1+v)^{2g} (1+ 2 u^{g+1} v^{g+1} -u^2v^2) -
  (1-u^2)^g(1-v^2)^g(1-uv)^2  \bigg)  \\ & +
  \frac12 \Big((1+u)^{2g}(1+v)^{2g}
   + (1-u^2)^{g}(1-v^2)^g \Big) \, .
 \end{align*}
\end{corollary}

\begin{proof}
  We only need to compute $e(M^{ss}(2,d))$, where
  $M^{ss}(2,d)=M(2,d)-M^s(2,d)$ is the locus of non-stable and
  polystable rank $2$ bundles of degree $d$. Such bundles are of
  the form $L_1\oplus L_2$, where $L_1,L_2\in \Jac^{d/2} X$.
  Therefore $M^{ss}(2,d) \cong (\Jac\ X\x \Jac\ X )/\ZZ_2$.
By Lemma \ref{lem:Z2} and (\ref{eqn:Jac}),
   $$
   e((\Jac\ X\x \Jac\ X )/\ZZ_2)= \frac12 \Big((1+u)^{2g}(1+v)^{2g}
   + (1-u^2)^{g}(1-v^2)^g \Big)\, .
   $$
Adding this to $e(M^s(2,d))$ in Theorem \ref{thm:rank2even} we get
the result.
\end{proof}

For instance, the formula of Corollary \ref{cor:rank2even-ss} for
$g=2$ gives
  $$
  e(M(2,0))=(1+u)^2(1+v)^2(1+uv+u^2v^2+u^3v^3)\, .
  $$
This formula agrees with \cite[Remark 4.11]{Ki}. Note that the
moduli space $M(2,0)$ is smooth for $g=2$ (see \cite{NR2}).

\section{Critical values for triples of rank $(2,2)$}\label{sec:critical(2,2)}

Now we move to the analysis of the moduli spaces of
$\s$-polystable triples of rank $(2,2)$. Let $\cN_\s=
\cN_\s(2,2,d_1,d_2)$. By Proposition \ref{prop:alpha-range}, $\s$
takes values in the interval
 $$
 I=[\s_m,\infty)=[\mu_1-\mu_2,\infty)\, , \quad \hbox{ where $d_1-d_2\geq 0$}.
 $$
Otherwise $\cN_\s$ is empty.

\begin{theorem} \label{thm:moduli(2,2)}
For $\s\in I$, $\cN_\s$ is a projective variety.  It is smooth of
dimension $4g+2d_1-2d_2-3$ at any $\s$-stable point for $\s\geq
2g-2$, or at any $\s$-stable injective triple. Moreover, if
$d_1+d_2$ is odd then $\cN_\s=\cN_\s^s$ for non-critical $\s$.
\end{theorem}

\begin{proof}
Projectiveness follows from Proposition \ref{prop:alpha-range}.
The smoothness at injective triples follows from Theorem
\ref{thm:smoothdim}(4); the dimension follows from Theorem
\ref{thm:smoothdim}(3); the smoothness result for $\s\geq 2g-2$
comes from \cite[Theorem 3.8(6)]{BGPG}. If $d_1+d_2$ is odd then
$\GCD(2,2,d_1+d_2)=1$ and so, for non-critical $\s$,
$\cN_\s=\cN_\s^s$, by Proposition
\ref{prop:triples-critical-range} (3). On the other hand, if
$d_1+d_2$ is even, then it may happen that there are strictly
$\s$-semistable triples for non-critical values of $\s$.
\end{proof}

Let us now compute the critical values for $\cN_\s(2,2,d_1,d_2)$.
According to (\ref{eqn:sigmac}) we have the following
possibilities for $n_1=2$, $n_2=2$:

\begin{enumerate}
 \item[(1)] $n'_1=1$, $n'_2=0$. The corresponding $\s_c$-destabilizing
 subtriple is of the form $0\to E_1'$ where $E_1'=L$ is a line
 bundle of degree $d_L$. The critical value is
   $$
   \s_c= \frac{4d_L-(d_1+d_2)}{2}=2d_L-\mu_1-\mu_2 \, .
   $$

\item[(2)] $n'_1=1$, $n'_2=2$. The $\s_c$-destabilizing subtriple
$T'$ is of the form $E_2 \to E_1'$ where $E_1'$ is a line bundle.
The quotient triple $T''=T/T'$ is of the form $0\to E_1''$, where
$E_1''=L$ is a line bundle of degree $d_L$, and $d_1'=d_1-d_L$.
Note that $\phi:E_2\to E_1$ is not injective. The critical value
is
 $$
 \s_c= \frac{4(d_1-d_L+d_2)-3(d_1+d_2)}{-2}= 2d_L-\mu_1-\mu_2 \, .
 $$

\item[(3)] $n'_1=2$, $n'_2=1$. The $\s_c$-destabilizing subtriple
$T'$
 is of the form $E_2'\to E_1$, where $E_2'$ is a line bundle. Then
 the quotient triple $T''=T/T'$ is of the form $E_2''\to 0$, where
 $E_2''=F$ is a line bundle of degree $d_F$, and $d_2'=d_2-d_F$.
  $$
  \s_c =\frac{4(d_1+d_2-d_F)-3(d_1+d_2)}{2}=\mu_1+\mu_2-2d_F \, .
  $$

 \item[(4)] $n'_1=0$, $n'_2=1$. The $\s_c$-destabilizing subtriple is
 of the form $E_2'\to 0$, where $E_2'=F$ is a line bundle of
 degree $d_F$. Again in this case $\phi$ is not injective. The
 corresponding critical value is
   $$
   \s_c= \frac{4d_F-(d_1+d_2)}{-2}=\mu_1+\mu_2-2d_F\, .
   $$

\item[(5)] $n'_1=2$, $n'_2=0$. The subtriple is of the form $0\to
E_1$. The corresponding critical value is $\s_c=\mu_1-\mu_2=\s_m$.

\item[(6)] $n'_1=0$, $n'_2=2$. The subtriple is of the form
$E_2\to 0$. This only happens if $\phi=0$, and so
$T=(0,E_2,0)\oplus (E_1,0,0)$. The critical value is
$\s_c=\mu_1-\mu_2=\s_m$, and the triple is $\s$-unstable for any
$\s\neq \s_m$.

\end{enumerate}

Note that the case $n'_1=1$, $n'_2=1$ does not appear, since
$\lambda'=\lambda$ and therefore this does not give a critical
value. In the Cases (1), (3) and (5), we have $\lambda'<\lambda$,
so the corresponding triples are $\s$-unstable for $\s<\s_c$. In
the Cases (2), (4) and (6), we have $\lambda'>\lambda$, so the
corresponding triples are $\s$-unstable for $\s>\s_c$.

\begin{proposition}\label{prop:bounds-for-sc}
 \begin{itemize}
  \item[(i)] Let $\s_c=2d_L-\mu_1-\mu_2$ be a critical value corresponding to the
  Cases {\rm (1)} or {\rm (3)}. Then $\mu_1 \leq d_L\leq (3\mu_1-\mu_2)/2$.
  Also $d_L =\mu_1 \iff \s_c=\s_m$.
  \item[(ii)] Let $\s_c=\mu_1+\mu_2-2d_F$ be a critical value corresponding to
  the Cases {\rm (2)} or {\rm (4)}. Then $(3\mu_2-\mu_1)/2 \leq d_F\leq
  \mu_2$. Also $d_F =\mu_2 \iff \s_c=\s_m$.
 \end{itemize}
\end{proposition}

\begin{proof}
We shall do the first item, since the second is analogous. Fix the
critical value $\s_c=2d_L-\mu_1-\mu_2$ and suppose that there is a
strictly $\s_c$-semistable triple $T$ in either Case (1) or (3)
above. Then the subtriple $T'$ and quotient triple $T''$ are both
$\s_c$-semistable by Lemma \ref{lem:semistable}. In either case,
there exists a $\s_c$-semistable triple of type
$(1,2,d_1-d_L,d_2)$. By Proposition \ref{prop:alpha-range} applied
to this situation, we get
 $$
    d_1-d_L-\frac{d_2}2 \leq \s_c= 2d_L-\frac{d_1}2-\frac{d_2}2 \leq
    4\left(d_1-d_L-\frac{d_2}2 \right)\, .
 $$
We can write this inequality in the equivalent form
 $$
    \frac{d_1}2 \leq d_L \leq \frac{3d_1-d_2}4 \,.
 $$
\end{proof}

\begin{theorem} \label{thm:stabilize}
  Let $\s_M= 2(\mu_1-\mu_2)$. For
  $\s>\s_M$ the moduli spaces of $\s$-(semi)stable triples do not change,
  and all $\s$-semistable triples
  $T=(E_1,E_2,\phi)$ are injective, i.e., $T$ defines an
  exact sequence of the form
  $$
   0 \to E_2 \overset{\phi}{\lto} E_1 \to S \to 0,
  $$
  where $S$ is a torsion sheaf of degree  $d_1 -d_2$.
\end{theorem}

\begin{proof}
  If we are in the first situation in Proposition
  \ref{prop:bounds-for-sc}, then $\s_c =2d_L-\mu_1-\mu_2 \leq
  3\mu_1-\mu_2-\mu_1-\mu_2=2(\mu_1-\mu_2)$. In the second
  situation, $\s_c=\mu_1+\mu_2-2d_F\leq \mu_1+\mu_2-
  (3\mu_2-\mu_1)=2(\mu_1-\mu_2)$.

  Now let $T$ be a $\s$-semistable triple for $\s>2(\mu_1-\mu_2)$.
  If $\phi:E_2\to E_1$ were not injective, then $T$ has a
  subtriple $T'=(0,\ker\phi,0)$ with $\lambda'>\lambda$. This
  forces $\mu_\s(T')>\mu_\s(T)$ for $\s$ large, and hence for $\s$
  bigger than the last critical value.
\end{proof}

\begin{remark}\label{rem:stabilize}
  Note that for any critical value $\s_c$, all the triples in
  $\Smenos$ are not injective.
\end{remark}

\begin{remark}\label{rem:stabilization}
 By \cite[Proposition 6.5]{BGPG} there is a value $\s_0$ such that
 all $\s$-semistable triples for $\s>\s_0$ are injective. By
 \cite[Theorem 8.6]{BGPG} there is a value $\s_L$ such that the
 moduli spaces $\cN_\s$ are isomorphic for all $\s>\s_L$.
 In our case, $n_1=n_2=2$,
 both numbers are $2(\mu_1-\mu_2)$.
\end{remark}

\begin{remark} \label{rem:virtual-not-all}
 In Proposition \ref{prop:bounds-for-sc} we see that, for the
 triples of rank $(2,2)$, there are critical values for which
 the moduli spaces do not change (those corresponding to
 $d_L>(3\mu_1-\mu_2)/2$ and those corresponding to
 $d_F<(3\mu_2-\mu_1)/2$).
\end{remark}

\begin{remark}\label{rem:simple-double}
  If we have simultaneously $\s_c=2d_L-\mu_1-\mu_2$ and
  $\s_c=\mu_1+\mu_2-2d_F$, then $2d_L-\mu_1-\mu_2=\mu_1+\mu_2-2d_F
  \implies d_1+d_2 =2d_L+2d_F$ is an even number.

  Therefore, if
  $d_1+d_2\notin 2\ZZ$, then Cases (1) and (3) (resp.\ Cases (2) and
  (4)) do not happen
  simultaneously (for the same critical value). So the flip locus $\Smas$
  (resp.\ $\Smenos$) will consist only of triples of one type for
  any $\s_c>\s_m$. In this situation
  the critical values $\s_c\in (\mu_1+\mu_2+ \ZZ) \cap
  [\mu_1-\mu_2, 2(\mu_1-\mu_2)]$.

  If $d_1+d_2\in 2\ZZ$, then Cases (1) and (3) (resp.\ Cases (2) and
  (4)) do happen simultaneously. The flip locus $\Smas$ (resp.\ $\Smenos$)
  consists of two types of triples, which yields two components
  that must be considered independently. In this situation
  the critical values $\s_c\in (\mu_1+\mu_2+ 2\ZZ) \cap
  [\mu_1-\mu_2, 2(\mu_1-\mu_2)]$.
\end{remark}

In the next section, it will be useful to have a vanishing result
for the hypercohomology $\HH^2$ to find the flip loci $\cS_\scpm$
for the moduli spaces of triples of type $(2,2,d_1,d_2)$.

\begin{proposition} \label{prop:h2-vanishing}
   Let $T=(E_1,E_2,\phi)$ be a strictly $\s_c$-semistable triple of
   type $(2,2,d_1,d_2)$ with $\s_c>\s_m$, $T'=(E_1',E_2',\phi')$ a destabilizing
   subtriple and $T''=T/T'=(E_1'',E_2'',\phi'')$ the corresponding
   quotient triple.
   \begin{itemize}
    \item[(1)] If $T \in \Smenos$ then $\HH^{2}(C^{\bullet}(T'',T')) = 0$.
    \item[(2)] If $T \in \Smas$ then $\HH^{2}(C^{\bullet}(T'',T')) = 0$,
     if $d_1-d_2 > 2g-2$.
   \end{itemize}
\end{proposition}

\begin{proof}
By Proposition \ref{prop:hyper-equals-hom} and Serre duality, the
vanishing $\HH^{2}(C^{\bullet}(T'',T'))=0$ is equivalent to the
injectivity of the map
 $$
 \begin{array}{ccc}
  H^{0}({E_{1}'}^* \otimes E_2'' \otimes K) & \overset{P}{\lto} &
  H^{0}({E_{1}'}^* \otimes E_1'' \otimes
  K)\oplus H^{0}({E_{2}'}^* \otimes E_2'' \otimes K) \\
 \psi & \longmapsto & ((\phi'' \otimes Id) \circ \psi, \, \psi
 \circ \phi').
 \end{array}
 $$

\begin{itemize}
\item[(1)]  If $T \in \Smenos$, then $H^{0}({E_{1}'}^* \otimes
E_2'' \otimes K)$ is trivial because either we are in Case (4) and
so $E_{1}'=0$ or we are in Case (2) and so $E_{2}''=0$.

\item[(2)] If $T \in \Smas$, we may have two cases:
\begin{itemize}
\item[(a)] If we are in Case (3), then $E_1'=E_1$ and $E_1''=0$.
The map $P$ is
 $$
 \begin{array}{ccc}
  H^{0}(E_{1}^* \otimes E_2'' \otimes K) & \overset{P}{\lto} &
  H^{0}({E_{2}'}^* \otimes E_2'' \otimes K) \\
 \psi & \longmapsto & \psi \circ \phi'.
 \end{array}
 $$
If $P$ is not injective, let $\psi: E_{1} \rightarrow E_2''
\otimes K$ be a non-trivial homomorphism in $\ker P$. Then, as
$\phi': E_2' \rightarrow E_1$, $\psi$ must factor through the
quotient $E_1/E_2'$. Both $E_1/E_2'$ and $E_2'' \otimes K$ are
line bundles, hence $\deg(E_1/E_2')=  d_1-d_2' \leq \deg(E_2''
\otimes K) =d_2''+2g-2$. This yields $d_1-d_2 \leq 2g-2$.

\item[(b)] If we are in Case (1), then $E_2'=0$ and $E_2''=E_2$.
Then the map $P$ is
 $$
 \begin{array}{ccc}
  H^{0}({E_{1}'}^* \otimes E_2'' \otimes K) & \overset{P}{\lto} &
  H^{0}({E_{1}'}^* \otimes E_1'' \otimes K) \\
  \psi & \longmapsto & (\phi'' \otimes Id) \circ \psi.
 \end{array}
 $$
If $P$ is not injective, let $\psi: E_{1}' \rightarrow E_2 \otimes
K$ be a non-trivial homomorphism in $\ker P$. Denote by $Q$ the
kernel of $\phi'': E_2 \rightarrow E_1''$, so $\psi$ must factor
through $Q \otimes K$. As $E_1'$ and $Q\otimes K$ are line
bundles, we have $\deg({E_1'})=d_1' \leq \deg(Q \otimes K) =
d_2-d_1''+2g-2$, which is rewritten as $d_1-d_2 \leq  2g-2$.

\end{itemize}
In both cases, if $P$ is not injective then $d_1-d_2 \leq 2g-2$.
Therefore, if $d_1-d_2> 2g-2$, then $P$ must be injective.
\end{itemize}
\end{proof}

\begin{remark}\label{rem:compareBGPG}
This result is a sort of improvement of \cite[Proposition
3.6]{BGPG} for the case of triples of rank $(2,2)$. Here we prove
the vanishing of $\HH^2$ for \textit{any} critical value $\s_c$
under the condition $\s_m=\mu_1-\mu_2> g-1$, whereas in
\cite[Proposition 3.6]{BGPG} it is proved the vanishing of $\HH^2$
only for critical values $\s_c>2g-2$ (but without condition in
$\s_m$).
\end{remark}

\section{Hodge polynomial of the moduli
of triples of rank $(2,2)$ and small $\s$} \label{sec:small}

In this section we want to compute the Hodge polynomial of the
moduli space
 $$
 \cN_{\smp}=\cN_{\smp}(2,2,d_1,d_2)
 $$
of $\s$-stable triples of types $(2,2,d_1,d_2)$ for $\s$ small,
under the assumption $\mu_1-\mu_2>2g-2$. The study of $\cN_{\smp}$
is simpler when both $d_1$ and $d_2$ are odd, since in this case
the bundles are automatically stable. However in this case
$d_1+d_2$ is even and hence $\gcd(2,2,d_1+d_2)\neq 1$. So there
may be strictly $\s$-semistable triples in $\cN_\s$ for
non-critical values of $\s$, making the moduli space $\cN_\s^s$
non-compact and the moduli space $\cN_\s$ singular (this does
not happen for $\s=\smp$; see Theorem \ref{thm:(2,2,odd,odd)}).

\begin{theorem}\label{thm:(2,2,odd,odd)}
Suppose that $d_1$ and $d_2$ are odd and that $\mu_1-\mu_2>2g-2$.
Then $\cN_\smp=\cN_\smp^s$, it is smooth, compact and
 $$
 e(\cN_\smp) = \left(\frac{(1+u)^{g}(1+v)^g(1+u^2v)^{g}(1+uv^2)^g
    -(uv)^{g}(1+u)^{2g}(1+v)^{2g}}
    {(1-uv)(1-(uv)^2)}\right)^2
    \frac{1-(uv)^{2d_1-2d_2-4g+4}}{1-uv} \, .
 $$
\end{theorem}

\begin{proof}
 The equality $\cN_\smp=\cN_\smp^s$ is a consequence of Proposition
 \ref{prop:moduli-small} (i).
 Next, since $\s_m=\mu_1-\mu_2>2g-2$, Theorem
 \ref{thm:moduli(2,2)} implies that
 the moduli $\modulimmas$ is smooth and compact.
 By Proposition \ref{prop:moduli-small} (i), it is
 the projectivization of a fiber bundle over $M(2,d_1) \times
 M(2,d_2)$ of rank $2d_1-2d_2-4g+4$. Therefore
 $$
 e(\cN_\smp) = e(M(2,d_1)) e(M(2,d_2)) e_{2d_1-2d_2-4g+4} \, .
 $$
 The result follows now applying Theorem \ref{thm:rank2odd}.
\end{proof}

The case where $d_1$ is odd and $d_2$ is even is more involved,
since we have to deal with the presence of strictly semistable
bundles in $M(2,d_2)$.

\begin{theorem}\label{thm:(2,2,odd,even)}
Suppose that $d_1$ is odd and $d_2$ is even and that
$\mu_1-\mu_2>2g-2$. Then $\cN_\smp=\cN_\smp^s$, it is smooth and
compact and
 $$
 \begin{aligned}
 e(\modulimmas) = & \ \frac{(1+u)^{2g}(1+v)^{2g}(1-(uv)^N)
 (u^gv^g(1+u)^g(1+v)^g-(1+u^2v)^g(1+uv^2)^g)}{(1-uv)^3(1-(uv)^2)^2}
 \cdot \\ &
 \bigg(
  (1+u)^{g} (1+v)^g (u^{g+1}v^{g+1} + u^{N+g-1}v^{N+g-1}) -
  (1+u^2v)^{g}(1+uv^2)^g (1 + u^Nv^{N})\bigg) \, ,
 \end{aligned}
  $$
where $N=d_1-d_2-2g+2$.
\end{theorem}

\begin{proof}
As $d_1+d_2$ is odd, Theorem \ref{thm:moduli(2,2)} implies that
$\modulimmas=\modulimmas^s$, and it is smooth and compact, since
$\s_m=\mu_1-\mu_2>2g-2$. To compute $e(\cN_\smp)$ we decompose
$\cN_\smp=X_0\sqcup X_1\sqcup X_2\sqcup X_3\sqcup X_4$, where:

\begin{enumerate}
\item[(1)] The open subset $X_0\subset \modulimmas$ consists of
those triples of the form $\phi:E_2\to E_1$, where $E_1$ and $E_2$
are both stable bundles, and $\phi$ is a non-zero map defined up
to multiplication by scalars. By Proposition
\ref{prop:moduli-small} (ii), $X_0 \to M(2,d_1)\times M^s(2,d_2)$
is a projective fibration whose fibers are projective spaces of
dimension $2d_1-2d_2-4g+4-1=2N-1$. Therefore, and using the
notation (\ref{eqn:Pn}),
  $$
  e(X_0)= e(M(2,d_1)) e(M^s(2,d_2))e_{2N} \, .
  $$

\item[(2)] The subset $X_1$ parametrizes $\smp$-stable triples of
the form $\phi:E_2\to E_1$ where $E_2$ is a strictly semistable
bundle of degree $d_2$ which is a non-split extension
 $$
  0 \to L_1 \to E_2\to L_2\to 0,
 $$
where $L_1, L_2\in \Jac^{d_2/2} X$ are non-isomorphic and $E_1$ is a
stable bundle. The space $Y_1$ parametrizing such bundles $E_2$ was
described in (2) of the proof of Theorem \ref{thm:rank2even} and its
Hodge polynomial is given in (\ref{eqn:Y1}).

Now in order to describe $X_1$, we must characterize when a triple
$T=(E_1,E_2,\phi)$, with $E_2\in Y_1$, is $\smp$-stable. As $T$ is
$\s_m$-semistable, then the only possibility for $T$ being
$\smp$-unstable is that it has a subtriple $T'$ of rank $(1,2)$ or
$(0,1)$, corresponding to Cases (2) or (4) of Section
\ref{sec:critical(2,2)}, respectively. If $T'$ is of rank $(1,2)$,
then it is of the form $E_2\to L$, where $L$ is a line bundle of
degree $d_L=\mu_1$, by Proposition \ref{prop:bounds-for-sc}. But
this is impossible, since $d_1$ is odd. If $T'$ is of rank $(0,1)$,
then it is of the form $F \to 0$, where $F$ is a line bundle of
degree $d_F=\mu_2$, by Proposition \ref{prop:bounds-for-sc}.
Therefore $F$ is a destabilizing subbundle for $E_2$. Since the only
destabilizing subbundle of $E_2$ is $L_1$, we have $F=L_1$. So it
must be $\phi(L_1)=0$. Any such $\phi$ lies in the image of the
inclusion $\Hom(L_2,E_1)\inc \Hom(E_2,E_1)$, under the natural
projection $E_2\to L_2$. This discussion implies that given
$(E_1,E_2)\in M(2,d_1)\times Y_1$, the morphisms $\phi$ giving rise
to $\smp$-stable triples $(E_1,E_2,\phi)$ are those in
 $$
 \Hom (E_2,E_1) - \Hom (L_2,E_1)\ .
 $$
Note that since the group of automorphisms of $E_1$ and $E_2$ are
both equal to $\CC^*$, $\phi$ is defined up to multiplication by
non-zero scalars. So the map $\pi:X_1\to M(2,d_1)\times Y_1$ is a
fibration with fiber over $(E_1,E_2)$ equal to
 \begin{equation}\label{eqn:333}
  \PP\Hom (E_2,E_1) - \PP \Hom (L_2,E_1)\ .
 \end{equation}
By Riemman-Roch, $\dim \Hom(E_2,E_1)=2d_1-2d_2-4g+4=2N$, since
$\mu_1-\mu_2>2g-2$ implies that $H^1(E_2^*\otimes E_1)=
H^0(E_1\otimes E_2^*\otimes K)=0$, $E_1$ and $E_2$ being both
semistable bundles. Also $\dim
\Hom(L_2,E_1)=d_1-2(d_2/2)-2g+2=d_1-d_2-2g+2=N$, since $\mu_1-\deg
L_2=\mu_1-d_2/2>2g-2$.  Hence (\ref{eqn:333}) is isomorphic to
$\PP^{2N-1}-\PP^{N-1}$. Therefore as in (2) of the proof of Theorem
\ref{thm:rank2even},
 \begin{align*}
 e(X_1)= & \ e(M(2,d_1)) e(Y_1) (e_{2N}-e_N) \\
  =& \ e(M(2,d_1))e(\Jac X)(e(\Jac X)-1)e_{g-1} (e_{2N}-e_{N}) \, .
 \end{align*}

\item[(3)] The subset $X_2$ parametrizes $\smp$-stable triples of
the form $\phi:E_2\to E_1$ where $E_2$ is a strictly semistable
bundle of degree $d_2$ which is non-split extension
 $$
  0 \to L_1 \to E_2\to L_1\to 0,
 $$
where $L_1 \in \Jac^{d_2/2} X$ and $E_1$ is a stable bundle. The
space $Y_2$ parametrizing such bundles $E_2$ was described in (3) of
the proof of Theorem \ref{thm:rank2even} and its Hodge polynomial is
given in (\ref{eqn:Y2}).

To describe $X_2$, we must characterize when a triple
$T=(E_1,E_2,\phi)$, with $E_2\in Y_2$, is $\smp$-stable. As before,
given $(E_1,E_2)\in M(2,d_1)\times Y_2$, the morphisms $\phi$ giving
rise to $\smp$-stable triples $(E_1,E_2,\phi)$ are those in
 $$
  \Hom (E_2,E_1) - \Hom (L_1,E_1)\ .
 $$

For a triple $T=(E_1,E_2,\phi)\in X_2$, $\Aut(E_1)=\CC^*$, so
$\Aut(T)\cong \Aut(E_2)=\CC\x\CC^*$. There is an exact sequence
 $$
 0\to \Hom(L_1,E_1) \to \Hom(E_2,E_1)\to \Hom(L_1,E_1)\to 0
 $$
Under the (non-canonical) decomposition $\Hom(E_2,E_1)\cong
\Hom(L_1,E_1)\oplus \Hom(L_1,E_1)$, $\Aut(E_2)$ acts as
$(a,\lambda) (x,y)\mapsto (\lambda x+a y, \lambda y)$. So the
fiber of $\pi:X_2\to M(2,d_1)\x Y_2$ is
 $$
  (\Hom (E_2,E_1) - \Hom (L_1,E_1))/\CC\x \CC^* \cong
  (\CC^{2N}-\CC^N)/\CC\x\CC^*\, ,
 $$
which is a $\CC^{N-1}$-bundle over $\PP^{N-1}$. Therefore
as in (3) of the proof of Theorem
\ref{thm:rank2even},
 \begin{align*}
 e(X_2)= & \ e(M(2,d_1)) e(Y_2)
 (e_{N}-e_{N-1})e_N \\ =& \ e(M(2,d_1)) e(\Jac X) e_g
 (e_{N}-e_{N-1})e_N \, .
 \end{align*}

\item[(4)] The subset $X_3$ parametrizes $\smp$-stable triples of
the form $\phi:E_2\to E_1$ where $E_1$ is a stable bundle and
$E_2=L_1\oplus L_2$, $L_1\not\cong L_2$, are two line bundles of
degree $d_2/2$. The space $Y_3$ parametrizing such bundles is
described in (\ref{eqn:Y3}).

As above, the condition for $\phi\in \Hom(E_2,E_1)$ to give rise to
a $\smp$-unstable triple is that there is a subtriple $T'$ of the
form $F\to 0$, with $F$ a line bundle of degree $d_F=\mu_2$. Then it
must be either $F=L_1$ or $F=L_2$. This means that $\phi\in
(\Hom(L_1,E_1)\oplus \{0\} )\cup (\{0\}\oplus\Hom(L_2,E_1)) \subset
\Hom(E_2,E_1)$. Therefore, given $(E_1,E_2)\in M(2,d_1)\times Y_3$,
the morphisms $\phi$ giving rise to $\smp$-stable triples
$(E_1,E_2,\phi)$ are those in
 $$
 (\Hom (L_1,E_1)-\{0\})\times (\Hom(L_2,E_1)-\{0\}).
 $$

The group of automorphisms of $E_2$ is $\CC^*\times \CC^*$ acting
on $L_1\oplus L_2$ by diagonal matrices. Therefore $\phi\in (\Hom
(L_1,E_1)-\{0\})\times (\Hom(L_2,E_1)-\{0\})$ is defined up to the
action of $\CC^*\times \CC^*$, where each $\CC^*$ acts by
multiplication on each of the two summands. So the map $\pi:X_3\to
M(2,d_1)\x Y_3$ has fiber
 \begin{equation}\label{eqn:444}
 \PP \Hom (L_1,E_1)\times\PP\Hom(L_2,E_1).
 \end{equation}
By Riemann-Roch, $\dim\Hom(L_1,E_1)=\dim\Hom(L_2,E_1)=
d_1-d_2-2g+2$. Therefore (\ref{eqn:444}) is isomorphic $\PP^{N-1}\x
\PP^{N-1}$. To compute $e(X_3)$ we work as in (4) of the proof of
Theorem \ref{thm:rank2even}. Write
$X_3=\tilde{X}_3/\ZZ_2=X_3'/\ZZ_2-X_3''/\ZZ_2$, where $X_3'$ is a
fibration over $M(2,d_1)$ with fiber $(A_{E_1}\x A_{E_1})/\ZZ_2$,
where $A_{E_1}$ is a projective bundle over $\Jac^{d_2/2}X$ with
fibers $\PP \Hom (L,E_1)\cong \PP^{N-1}$, and $\ZZ_2$ acts by
permutation. $X_3''$ is a fibration over $M(2,d_1)\x \Jac^{d_2/2}X$
with fibers $(\PP^{N-1}\x \PP^{N-1})/\ZZ_2$. So using Theorem
\ref{thm:Du},
  $$
  \begin{aligned}
  e(X_3)= & \ e(\tilde{X}_3/\ZZ_2)=  e(X_3'/\ZZ_2)-e(X_3''/\ZZ_2) \\
   = & \ \frac12 e(M(2,d_1)\bigg(\bigg( e(\Jac X)^2e_N^2 +  (1-u^2)^g(1-v^2)^g
   \frac{1-(uv)^{2N}}{1-u^2v^2}\bigg)  \\ &
   -e(\Jac X) \bigg( e_N^2 +
   \frac{1-(uv)^{2N}}{1-u^2v^2}\bigg)\bigg)\, .
  \end{aligned}
  $$

\item[(5)] The subset $X_4$ parametrizes triples $\phi:E_2\to
E_1$, where $E_1$ is a stable bundle and $E_2=L_1\oplus L_1$,
$L_1\in \Jac^{d_2/2}X$. Such bundles $E_2$ are parametrized by
$Y_4=\Jac^{d_2/2} X$. The map $\phi$ lies in
 \begin{equation}\label{eqn:555}
 \Hom(E_2,E_1)=\Hom(L_1,E_1)\oplus \Hom(L_1,E_1) \cong
 \Hom(L_1,E_1)\otimes \CC^2\, .
 \end{equation}

The condition for a triple $T=(E_1,E_2,\phi)$ to be $\smp$-unstable
is that there is a line subbundle $F\subset E_2$ of degree
$d_F=\mu_2$ such that $\phi(F)=0$. A destabilizing subbundle of
$E_2$ is necessarily isomorphic to $L$ and there exists $(a,b)\neq
(0,0)$ such that $F\cong L\inc E_2$ is given by $x\mapsto (ax,bx)$.
So $\phi=(a\psi,b\psi)\in\Hom(L_1,E_1)\otimes \CC^2$, for some
$\psi\in \Hom(L,E_1)$. Therefore $T=(E_1,E_2,\phi)$ is $\smp$-stable
if $\phi=(\phi_1,\phi_2)\in\Hom(L_1,E_1)\otimes \CC^2$ satisfies
that $\phi_1,\phi_2$ are linearly independent.

On the other hand, a triple $(E_1,E_2,\phi)\in X_4$ is determined up
to the action of $\Aut(E_2)=GL(2,\CC)$. This acts on (\ref{eqn:555})
via the standard representation on $\CC^2$. Thus for $(E_1,E_2)\in
M(2,d_1)\times Y_4$, the morphisms $\phi$ giving rise to
$\smp$-stable triples $(E_1,E_2,\phi)$ are parametrized by
$\Gr(2,\Hom(L_1,E_1))$. But $\dim \Hom(L_1,E_1)=d_1-d_2-2g+2=N$, so
this fiber is isomorphic to $\Gr(2,N)$. So
  $$
  e(X_4)=e(M(2,d_1)) e(Y_4) e(\Gr(2,N))=e(M(2,d_1)) e(\Jac X) e(\Gr(2,N)) \ .
  $$

\end{enumerate}

Adding up all contributions together we get
  \begin{align*}
  e(\cN_{\smp})= & \ e(X_0)+ e(X_1)+e(X_2)+e(X_3)+e(X_4) \\
  =& \ e(M(2,d_1)) \Bigg( e(M^s(2,d_2))e_{2N} + e(\Jac X) (e(\Jac X)-1) e_{g-1}
  (e_{2N}-e_{N})
  \\
  &\qquad+ e(\Jac X) e_g (e_{N}-e_{N-1})e_N  
  +  \frac12 \bigg( e(\Jac X)^2e_N^2 +  (1-u^2)^g(1-v^2)^g
   \frac{1-(uv)^{2N}}{1-u^2v^2}\bigg)
   \\
   &\qquad
   -\frac12 e(\Jac X) \bigg( e_N^2 +
   \frac{1-(uv)^{2N}}{1-u^2v^2}\bigg)    
  + e(\Jac X) e(\Gr(2,N)) \Bigg) \\
  =& \ \frac{(1+u)^{2g}(1+v)^{2g}(1-(uv)^N)
 (u^gv^g(1+u)^g(1+v)^g-(1+u^2v)^g(1+uv^2)^g)}{(1-uv)^3(1-(uv)^2)^2} \cdot
 \\ &
 \bigg(
  (1+u)^{g} (1+v)^g (u^{g+1}v^{g+1} +  u^{N+g-1}v^{N+g-1}) -
  (1+u^2v)^{g}(1+uv^2)^g (1 + u^Nv^{N})\bigg)  \, .
  \end{align*}
\end{proof}

\begin{corollary}\label{cor:(2,2,odd,even)-2}
Suppose that $d_1$ is even and $d_2$ is odd and that
$\mu_1-\mu_2>2g-2$. Then $\modulimmas=\modulimmas^s$, it is smooth
and compact and its Hodge polynomial has the same formula as that
of Theorem \ref{thm:(2,2,odd,even)}, where $N=d_1-d_2-2g+2$.
\end{corollary}

\begin{proof}
We use the isomorphism $\cN_{\s}(2,2,d_1,d_2)\cong
\cN_{\s}(2,2,-d_2,-d_1)$. Note that
  $$
  d_1-d_2= (-d_2) - (-d_1)\, ,
  $$
so that the small value $\smp=\mu_1-\mu_2$ and the condition on
the slopes $\mu_1-\mu_2>2g-2$ is the same for both moduli spaces
$\cN_{\s}(2,2,d_1,d_2)$ and $\cN_{\s}(2,2,-d_2,-d_1)$. Now we
apply Theorem \ref{thm:(2,2,odd,even)} to get the stated formula
where $N=-d_2-(-d_1)-2g+2$.
\end{proof}

\begin{corollary}\label{cor:(2,2,odd,even)}
Suppose that $d_1+d_2$ is odd and
$\mu_1-\mu_2>2g-2$. Then the Poincar{\'e} polynomial of $\cN_\smp$ is
 $$
 P_t(\modulimmas) = \frac{(1\!+\!t)^{4g} (1\!-\!t^{2N})(t^{2g}(1\!+\!t)^{2g} \!-\! (1\!+\!t^3)^{2g})
  ((1\!+\!t)^{2g} (t^{2g\!+\!2}\!+\! t^{2N\!+\!2g\!-\!2}) \!-\! (1\!+\!t^3)^{2g}
  (1 \!+\! t^{2N}))}{(1\!-\!t^2)^3(1\!-\!t^4)^2} \, ,
 $$
where $N=d_1-d_2-2g+2$. \hfill $\Box$
\end{corollary}

\begin{proof}
 $\cN_{\smp}$ is smooth and projective, so $P_t(\cN_{\smp})=e(\cN_{\smp})(t,t)$.
 The result follows from Theorem \ref{thm:(2,2,odd,even)} and
 Corollary \ref{cor:(2,2,odd,even)-2}.
\end{proof}

We could deal also with the case when $d_1$ and $d_2$ are both even
and $d_1-d_2>4g-4$. This is similar to the case just treated in
Theorem \ref{thm:(2,2,odd,even)}, with the further complication that
there are semistable loci for both $E_1$ and $E_2$.

However, dealing with the case $d_1-d_2\leq 4g-4$ is more
complicated, since Proposition \ref{prop:moduli-small} does not
apply as there is a Brill-Noether problem consisting on
determining the loci of those $(E_1,E_2)$ where $\dim
\Hom(E_2,E_1\ox K)$ is constant.

\section{Contribution of the flips to the Hodge polynomials} \label{sec:simple}

In this section, we shall compute the change in the Hodge
polynomial of $\cN_\s(2,2,d_1,d_2)$ when we cross a critical value
$\s_c$. We restrict to the case $d_1+d_2$ is odd, since in the
case $d_1+d_2$ even there may be strictly $\s$-semistable triples
for non-critical values of $\s$ (and in this case $\cN_\s^s$ is
non-compact and $\cN_\s$ is non-smooth). For $d_1+d_2$ odd,
Theorem \ref{thm:moduli(2,2)} guarantees that $\cN_\s$ is compact
and smooth for any non-critical $\s\geq 2g-2$.
The critical values are given in Proposition \ref{prop:bounds-for-sc}. These
are of two types. The following two propositions treat them separately.

\begin{proposition}\label{prop:diff-poly-simple-2}
Let $\s_c=2d_L-\mu_1-\mu_2$ be a critical value for triples of
type $(2,2,d_1,d_2)$ with $d_1+d_2$ odd, such that $\s_c>\s_m$.
Suppose that $\mu_1-\mu_2>g-1$. Then
 $$
 \begin{aligned}
  e(\modulimas)-e(\modulimenos)=
    \coeff_{x^0}\Bigg[&
    \frac{(1+u)^{3g}(1+v)^{3g}(1+ux)^{g}(1+vx)^g\big((uv)^{g-1-d_1+2d_L}-(uv)^{1-g+d_1-d_2}\big)}
    {(1-uv)^2(1-x)(1-uvx)x^{[3\mu_1-\mu_2]-2d_L}}\\
    &
    \Bigg(\frac{(uv)^{(3d_1-d_2-1)/2-2d_L}}{1-(uv)^{-1}x}
    -\, \frac{(uv)^{2d_L-d_1+g}}{1-(uv)^2x}\Bigg) \Bigg]\,.
 \end{aligned}
 $$
\end{proposition}

\begin{proof}
Theorem \ref{thm:moduli(2,2)} implies that
$\cN_\scpm=\cN_\scpm^s$. Then Lemma \ref{lem:fliploci} and the
properties of the Hodge polynomials give
 $$
  e(\cN_{\scp})-e( \cN_{\scm}) = e(\cS_{\scp}) -e(\cS_{\scm}).
 $$

Let us start by studying $\Smas$. By Lemma \ref{lem:semistable},
any $T\in \Smas$ sits in a non-split extension
 \begin{equation}\label{eqn:sec8}
    0 \rightarrow T' \rightarrow T \rightarrow T'' \rightarrow 0
 \end{equation}
in which $T'$ and $T''$ are $\s_c$-semistable, $\lambda'<\lambda$
and $\mu_{\s_c}(T')=\mu_{\s_c}(T)=\mu_{\s_c}(T'')$. Since $T$
corresponds to Case (1) in Section \ref{sec:critical(2,2)}, we
have $T'\in \cN_{\s_c}'$ and $T''\in \cN_{\s_c}''$, where
 \begin{align*}
  \cN_{\s_c}' &= \cN_{\s_c}(1,0,d_L,0) \cong \Jac^{d_L} X, \\
  \cN_{\s_c}''&= \cN_{\s_c}(1,2,d_1-d_L,d_2).
 \end{align*}
The moduli space of triples of rank $(1,0)$ has no critical
values; and for the moduli space of triples of rank $(1,2)$, the
critical values are of the form $3d_M + d_1''+d_2''$, by Lemma
\ref{lem:dM-2}, and are in particular integers. But
$\s_c=2d_L-\frac{d_1+d_2}{2} \notin\ZZ$, so $\s_c$ is not a
critical value for $\cN_{\s_c}''$.

By \cite[Proposition 3.5]{BGPG}, $\HH^0(T'',T')=0$ and
by Proposition \ref{prop:h2-vanishing} (2), $\HH^2(T'',T')=0$ . So
Theorem \ref{thm:Smas} implies that $\Smas$ is the
projectivization of a bundle over $\cN_{\s_c}'\times \cN_{\s_c}''$
of rank
 $$
 -\chi(T'',T')= 1 -g+d_1-d_2\, .
 $$
Therefore
 $$
 e(\Smas)= e(\Jac^{d_L} X)\, e(\cN_{\s_c}(1,2,d_1-d_L,d_2))\, e_{1 -g+d_1-d_2}\, .
 $$

The case of $\Smenos$ is similar. Any $T\in\Smenos$ sits in an
exact sequence (\ref{eqn:sec8}) with $T'\in \cN_{\s_c}'$ and
$T''\in \cN_{\s_c}''$, where
 \begin{align*}
  \cN_{\s_c}' &= \cN_{\s_c}(1,2,d_1-d_L,d_2), \\
  \cN_{\s_c}''&= \cN_{\s_c}(1,0,d_L,0)\cong \Jac^{d_L} X,
 \end{align*}
corresponding to the Case (2) in Section \ref{sec:critical(2,2)}.
The hypothesis of Theorem \ref{thm:Smas} are satisfied and so
$\Smenos$ is the projectivization of a bundle over
$\cN_{\s_c}'\times \cN_{\s_c}''$ of rank
 $$
 -\chi(T'',T')= g-1 -d_1+2d_L\, .
 $$
Therefore
 $$
 e(\Smenos)= e(\Jac^{d_L} X) \, e(\cN_{\s_c}(1,2,d_1-d_L,d_2)) \, e_{g-1 -d_1+2d_L}\, .
 $$

Substracting, we get
 \begin{align*}
 e(\Smas)-e(\Smenos) &= (e_{1
 -g+d_1-d_2}- e_{g-1 -d_1+2d_L}) (1+u)^{g}(1+v)^g
 e(\cN_{\s_c}(1,2,d_1-d_L,d_2)) = \\
  &= \frac{(uv)^{g-1-d_1+2d_L}-(uv)^{1-g+d_1-d_2}}{1-uv}
  (1+u)^{g}(1+v)^g e(\cN_{\s_c}(1,2,d_1-d_L,d_2))\, .
 \end{align*}

Being $\s_c$ a non-critical value for the moduli of triples of
rank $(1,2)$, we can apply Theorem
\ref{thm:polinomiono(1,2)no-critico} to compute the Hodge
polynomial of $\cN_\s(1,2,d_1-d_L,d_2)$. First,
  $$
  \begin{aligned}
  d_0 &=\left[\frac13 (2d_L-\mu_1-\mu_2-(d_1-d_L)-d_2)\right]+1\\
   &=d_L +[-\mu_1-\mu_2] +1 \, .
  \end{aligned}
  $$
So $e(\cN_\s(1,2,d_1-d_L,d_2))$ equals
 $$
 \coeff_{x^0} \Bigg[
    \frac{(1+u)^{2g}(1+v)^{2g}(1+ux)^{g}(1+vx)^{g}}{(1-uv)(1-x)(1-uv x)x^{d_1-d_2-d_L-d_0}}
    \Bigg(\frac{(uv)^{d_1-d_2-d_L-d_0}}{1-(uv)^{-1}x}-\,
    \frac{(uv)^{d_2+g-1+2d_0}}{1-(uv)^2x}\Bigg) \Bigg]\, ,
 $$
where $d_1-d_2-d_L-d_0= [3\mu_1-\mu_2]-2d_L=(3d_1-d_2-1)/2-2d_L$
and $d_2+2d_0=2d_L-d_1+1$. The result follows from this.
\end{proof}

\begin{proposition}\label{prop:diff-poly-simple-1}
Let $\s_c=\mu_1+\mu_2 - 2d_F$ be a critical value for triples of
type $(2,2,d_1,d_2)$ with $d_1+d_2$ odd, such that $\s_c>\s_m$.
Suppose that $\mu_1-\mu_2>g-1$. Then
 $$
 \begin{aligned}
  e(\modulimas)-e(\modulimenos)=
    \coeff_{x^0}\Bigg[&
    \frac{(1+u)^{3g}(1+v)^{3g}(1+ux)^{g}(1+vx)^g\big((uv)^{g-1+d_2-2d_F}-(uv)^{1-g+d_1-d_2}\big)}
    {(1-uv)^2(1-x)(1-uvx)x^{2d_F - [3\mu_2-\mu_1]-1}}\\
    &
    \Bigg(\frac{(uv)^{2d_F+(d_1-3d_2-1)/2}}{1-(uv)^{-1}x}
    -\, \frac{(uv)^{d_2-2d_F+g}}{1-(uv)^2x}\Bigg) \Bigg] \,.
 \end{aligned}
 $$
\end{proposition}

\begin{proof}
This is very similar to the proof of Proposition
\ref{prop:diff-poly-simple-2}. Again
 $$
  e(\cN_{\scp})-e( \cN_{\scm}) = e(\cS_{\scp})
  -e(\cS_{\scm}).
 $$
We start with $\Smas$. Any $T\in \Smas$ sits in a non-split
extension like (\ref{eqn:sec8}), with
$\mu_{\s_c}(T')=\mu_{\s_c}(T)=\mu_{\s_c}(T'')$, $T'\in
\cN_{\s_c}'$ and $T''\in \cN_{\s_c}''$, where
 \begin{align*}
  \cN_{\s_c}' &= \cN_{\s_c}(2,1,d_1,d_2-d_F), \\
  \cN_{\s_c}''&= \cN_{\s_c}(0,1,0,d_F)\cong \Jac^{d_F} X,
 \end{align*}
corresponding to the Case (3) in Section \ref{sec:critical(2,2)}.
The moduli space of triples of rank $(0,1)$ has no critical
values; and for the moduli space of triples of rank $(2,1)$, the
critical values are of the form $3d_M -d_1' -d_2' \in\ZZ$, whilst
$\s_c=\frac{d_1+d_2}{2} - 2d_F \notin\ZZ$, so $\s_c$ is not a
critical value for $\cN_{\s_c}'$. The other conditions of Theorem
\ref{thm:Smas} are checked as before. So $\Smas$ is the
projectivization of a bundle over $\cN_{\s_c}'\times \cN_{\s_c}''$
of rank
 $$
 -\chi(T'',T')= 1 -g+d_1-d_2\, .
 $$
Therefore
 $$
 e(\Smas)= e(\Jac^{d_F} X)\, e(\cN_{\s_c}(2,1,d_1,d_2-d_F))\, e_{1 -g+d_1-d_2}\, .
 $$

Moving to $\Smenos$, any $T\in\Smenos$ sits in an exact sequence
(\ref{eqn:sec8}) with $T'\in \cN_{\s_c}'$ and $T''\in
\cN_{\s_c}''$, where
 \begin{align*}
  \cN_{\s_c}' &= \cN_{\s_c}(0,1,0,d_F)\cong \Jac^{d_F} X, \\
  \cN_{\s_c}''&= \cN_{\s_c}(2,1,d_1,d_2-d_F),
 \end{align*}
corresponding to the Case (4) in Section \ref{sec:critical(2,2)}.
Arguing as before, we have that $\Smenos$ is the projectivization
of a bundle over $\cN_{\s_c}'\times \cN_{\s_c}''$ of rank
 $$
 -\chi(T'',T')= g-1 +d_2-2d_F\, .
 $$
Therefore
 $$
 e(\Smenos)= e(\Jac^{d_F} X) \, e(\cN_{\s_c}(2,1,d_1,d_2-d_F)) \,
 e_{g-1 +d_2-2d_F}\, .
 $$

Substracting, we get
 \begin{align*}
 e(\Smas)-e(\Smenos) &= (e_{1
 -g+d_1-d_2}-e_{g-1 +d_2-2d_F}) (1+u)^{g}(1+v)^g
 e(\cN_{\s_c}(2,1,d_1,d_2-d_F)) = \\
  &= \frac{(uv)^{g-1+d_2-2d_F}-t^{1-g+d_1-d_2}}{1-uv}
  (1+u)^{g}(1+v)^g e(\cN_{\s_c}(2,1,d_1,d_2-d_F))\, .
 \end{align*}

Being $\s_c$ a non-critical value for the moduli of triples of rank
$(2,1)$, we can apply Theorem \ref{thm:polinomiono(2,1)no-critico}
to compute the Hodge polynomial of $\cN_\s(2,1,d_1,d_2-d_F)$. First,
  $$
  \begin{aligned}
  d_0 &=\left[\frac13 (\mu_1+\mu_2-2d_F
  +d_1+d_2-d_F\right]+1\\
   &=[\mu_1+\mu_2] -d_F +1
  \end{aligned}
  $$
So $e(\cN_\s(2,1,d_1,d_2-d_F))$ equals
  $$
  \coeff_{x^0} \Bigg[
    \frac{(1+u)^{2g}(1+v)^{2g}(1+ux)^{g}(1+vx)^{g}}{(1-uv)(1-x)(1-uv x)x^{d_1-d_2+d_F-d_0}}
   \Bigg(\frac{(uv)^{d_1-d_2+d_F-d_0}}{1-(uv)^{-1}x}-\,
    \frac{(uv)^{-d_1+g-1+2d_0}}{1-(uv)^2x}\Bigg) \Bigg]\, ,
  $$
where $d_1-d_2+d_F-d_0=2d_F - [3\mu_2-\mu_1]-1=
2d_F+(d_1-3d_2-1)/2$ and $-d_1+2d_0=d_2-2d_F+1$. The result
follows from this.
\end{proof}

We gather together Propositions \ref{prop:diff-poly-simple-2} and
\ref{prop:diff-poly-simple-1} in a single result.

\begin{corollary}\label{cor:diff-poly-simple}
 The critical values $\s_c>\s_m$ for triples of type $(2,2,d_1,d_2)$ with
 $d_1+d_2$ odd are of the form $\s_c=\mu_1-\mu_2 +n$,
 $1\leq n \leq [\mu_1-\mu_2]$, $n\in \ZZ$. Suppose that $\mu_1-\mu_2>g-1$.
 Then
 $$
 \begin{aligned}
  e(\modulimas)-e(\modulimenos)=
    \coeff_{x^0}\Bigg[&
    \frac{(1+u)^{3g}(1+v)^{3g}(1+ux)^{g}(1+vx)^g\big((uv)^{g-1+n}-(uv)^{1-g+d_1-d_2}\big)}
    {(1-uv)^2(1-x)(1-uvx)x^{[\mu_1-\mu_2]-n}}\\
    &
    \Bigg(\frac{(uv)^{(d_1-d_2-1)/2-n}}{1-(uv)^{-1}x}
    -\, \frac{(uv)^{g+n}}{1-(uv)^2x}\Bigg) \Bigg]\,.
 \end{aligned}
 $$
\end{corollary}

\begin{proof}
 For simplicity let us assume that $d_1$ is odd and $d_2$ is even
 (the other case is analogous). We~have the
 following possibilities:
 \begin{enumerate}
 \item[(a)] If
 $\s_c=2d_L-\mu_1-\mu_2$, write $d_L =\mu_1+\frac12+m$ with $m$
 integer. Then $\s_c=\mu_1-\mu_2+2m+1$. As~$\mu_1<d_L\leq
 \frac{3\mu_1-\mu_2}2$ by Proposition \ref{prop:bounds-for-sc} (i),
 we have $0\leq m\leq (\mu_1-\mu_2-1)/2$. Substituting the
 values $3d_1-d_2-1-4d_L= d_1-d_2-1-4m-2$, $2d_L-d_1+g= g+2m+1$,
 $[3\mu_1-\mu_2]-2d_L= [\mu_1-\mu_2]-2m-1$ and $g-1-d_1+2d_L= g+2m$
 into the formula of Proposition \ref{prop:diff-poly-simple-2}, one
 gets
 $$
 \begin{aligned}
  e(\modulimas)-e(\modulimenos)=
    \coeff_{x^0}\Bigg[&
    \frac{(1+u)^{3g}(1+v)^{3g}(1+ux)^{g}(1+vx)^g\big((uv)^{g+2m}-(uv)^{1-g+d_1-d_2}\big)}
    {(1-uv)^2(1-x)(1-uvx)(1-(uv)^{-1}x)x^{[\mu_1-\mu_2]-2m-1}}\\
    &
    \Bigg(\frac{(uv)^{(d_1-d_2-1)/2-2m-1}}{1-(uv)^{-1}x}
    -\, \frac{(uv)^{g+2m+1}}{1-(uv)^2x}\Bigg) \Bigg]\,.
 \end{aligned}
 $$

\item[(b)] If $\s_c=\mu_1+\mu_2-2d_F$, write $d_F =\mu_2-m-1$ with $m$ an integer.
Then $\s_c=\mu_1-\mu_2+2m+2$. As~$\frac{3\mu_2-\mu_1}2 \leq d_F< \mu_2$  by
Proposition \ref{prop:bounds-for-sc} (i), we have $0\leq m\leq (\mu_1-\mu_2)/2-1$.
Substituting the values $4d_F+d_1-3d_2-1=d_1-d_2-1-4m-4$, $d_2-2d_F+g=g+2m+2$, $2d_F
- [3\mu_2-\mu_1]-1=[\mu_1-\mu_2] -2m-2$ and $g-1+d_2-2d_F= g+2m+1$ into the formula
of Proposition \ref{prop:diff-poly-simple-1}, we~have
 $$
 \begin{aligned}
  e(\modulimas)-e(\modulimenos)=
    \coeff_{x^0}\Bigg[&
    \frac{(1+u)^{3g}(1+v)^{3g}(1+ux)^{g}(1+vx)^g\big((uv)^{g+2m+1}-(uv)^{1-g+d_1-d_2}\big)}
    {(1-uv)^2(1-x)(1-uvx)(1-(uv)^{-1}x)x^{[\mu_1-\mu_2]-2m-2}}\\
    &
    \Bigg(\frac{(uv)^{(d_1-d_2-1)/2-2m-2}}{1-(uv)^{-1}x}
    -\, \frac{(uv)^{g+2m+2}}{1-(uv)^2x}\Bigg) \Bigg]\,.
 \end{aligned}
 $$
\end{enumerate}
Case (a) corresponds to $n=2m+1$ odd, and Case (b) to $n=2m+2$
even in the formula in the statement. The range for $n$ is $1\leq
n\leq \mu_1-\mu_2$. But, since $\mu_1-\mu_2$ is not an integer,
this range is actually $1\leq n \leq [\mu_1-\mu_2]$.
\end{proof}

\section{Hodge polynomial of the moduli
of triples of rank $(2,2)$ and large $\s$} \label{sec:large}

Now we use all the information in Sections
\ref{sec:critical(2,2)}--\ref{sec:simple}\  to compute the Hodge
polynomial of the $\moduli(2,2,d_1,d_2)$, for any non-critical
$\s>\s_m$. Recall that by Theorem \ref{thm:stabilize}, there is a
value $\s_M= 2(\mu_1-\mu_2)$ such that for $\s>\s_M$ all the moduli
spaces $\cN_\s$ are isomorphic. We refer to
 $$
 \cN_{\s_M^+}=\cN_{\s_M^+}(2,2,d_1,d_2)
 $$
as the \emph{large $\s$} moduli space.

\begin{proposition}\label{prop:finally}
 Suppose that $d_1$ is even and $d_2$ is odd and that $\mu_1-\mu_2>g-1$.
 Let $\s>\s_m$ be a non-critical value. Set
 $n_0=\min\{[\s-\mu_1+\mu_2], [\mu_1-\mu_2]\}$. Then
  $$
 \begin{aligned}
   e(\cN_{\s}) & - e(\cN_{\smp}) =   \coeff_{x^0}\Bigg[
    \frac{(1+u)^{3g}(1+v)^{3g}(1+ux)^{g}(1+vx)^g}
    {(1-uv)^2(1-x)(1-uvx)x^{[\mu_1-\mu_2]}}\\
    & \quad
     \Bigg(\frac{(uv)^{g-1+(d_1-d_2-1)/2}x(1-x^{n_0})}{(1-(uv)^{-1}x)(1-x)}
     -\, \frac{(uv)^{(3d_1-3d_2-1)/2-g}x (1- (uv)^{-n_0}
     x^{n_0})}{(1-(uv)^{-1}x)^2}\\
    &\quad -\, \frac{(uv)^{2g+1} x (1-(uv)^{2n_0}x^{n_0})}{(1-(uv)^2x)^2}
    +\, \frac{(uv)^{d_1-d_2+2}x (1-(uv)^{n_0} x^{n_0})}{(1-(uv)^2x)(1-uvx)}\Bigg) \Bigg]  \, .
  \end{aligned}
   $$
\end{proposition}

\begin{proof}
By Corollary \ref{cor:diff-poly-simple}, the critical values are
of the form $\s_c=\mu_1-\mu_2 +n$ with $1\leq n \leq
[\mu_1-\mu_2]$. Now $\s_m<\s_c<\s$ is equivalent to $n\leq
[\s-\mu_1+\mu_2]$ (note that $\s-\mu_1+\mu_2\not\in \ZZ$ since
$\s$ is not critical). Therefore,
 $$
 \begin{aligned}
 e(\cN_{\s})- & e(\cN_{\smp}) = \sum_{\s_m<\s_c<\s}
  e(\modulimas)-e(\modulimenos)= \\ = &\sum_{n=1}^{n_0}
    \coeff_{x^0}\Bigg[
    \frac{(1+u)^{3g}(1+v)^{3g}(1+ux)^{g}(1+vx)^g\big((uv)^{g-1+n}-(uv)^{1-g+d_1-d_2}\big)}
    {(1-uv)^2(1-x)(1-uvx)x^{[\mu_1-\mu_2]-n}}\\
    & \quad
    \Bigg(\frac{(uv)^{(d_1-d_2-1)/2-n}}{1-(uv)^{-1}x}
    -\, \frac{(uv)^{g+n}}{1-(uv)^2x}\Bigg) \Bigg] = \\
&= \coeff_{x^0}\Bigg[
    \frac{(1+u)^{3g}(1+v)^{3g}(1+ux)^{g}(1+vx)^g}
    {(1-uv)^2(1-x)(1-uvx)x^{[\mu_1-\mu_2]}}\\
    & \quad
     \Bigg(\frac{1}{1-(uv)^{-1}x} \sum_{n=1}^{n_0}(uv)^{g-1+(d_1-d_2-1)/2}x^n-
     \frac{1}{1-(uv)^{-1}x} \sum_{n=1}^{n_0}(uv)^{1-g+(3d_1-3d_2-1)/2-n}
     x^n \\ & \quad
    -\, \frac{1}{1-(uv)^2x} \sum_{n=1}^{n_0}(uv)^{2g-1+2n}x^n +
    \frac{1}{1-(uv)^2x} \sum_{n=1}^{n_0}(uv)^{1+d_1-d_2+n} x^n\Bigg) \Bigg] =
    \end{aligned}
 $$
 $$
 \begin{aligned}
   \qquad    &=
    \coeff_{x^0}\Bigg[
    \frac{(1+u)^{3g}(1+v)^{3g}(1+ux)^{g}(1+vx)^g}
    {(1-uv)^2(1-x)(1-uvx)x^{[\mu_1-\mu_2]}}\\
    & \quad
     \Bigg(\frac{(uv)^{g-1+(d_1-d_2-1)/2}x(1-x^{n_0})}{(1-(uv)^{-1}x)(1-x)}
     -\, \frac{(uv)^{1-g+(3d_1-3d_2-1)/2-1}x (1- (uv)^{-n_0}
     x^{n_0})}{(1-(uv)^{-1}x)^2}\\
    &\quad -\, \frac{(uv)^{2g-1+2} x (1-(uv)^{2n_0}x^{n_0})}{(1-(uv)^2x)^2}
    +\, \frac{(uv)^{1+d_1-d_2+1}x (1-(uv)^{n_0} x^{n_0})}{(1-(uv)^2x)(1-uvx)}\Bigg) \Bigg]\ .
 \end{aligned}
 $$
\end{proof}

\begin{theorem} \label{thm:finally}
 Suppose that $d_1$ is odd and $d_2$ is even.
 Then the large $\s$ moduli space $\cN_{\s_M^+}=\cN^s_{\s_M^+}$ is smooth
 and compact. If $\mu_1-\mu_2>2g-2$, its Hodge polynomial is
 $$
 \begin{aligned}
 e(\cN_{\s_M^+})=& \ \frac{(1+u)^{2g}(1+v)^{2g}}{(1-uv)^3(1-(uv)^2)^2}
    \Bigg[
    (1+u^2v)^{2g}(1+uv^2)^{2g}(1-(uv)^{2N})\\
    &
    -N \, (1+u^2v)^g(1+uv^2)^g(1+u)^g(1+v)^g (uv)^{N+g-1}(1-(uv)^2)
    \\
    &
    +
    (1+u)^{2g}(1+v)^{2g}(1+uv)^2(uv)^{2g -2 +(N+1)/2 } 
    \Big((1-(uv)^{N+1}) - \frac{N+1}{2} \, (1-uv)(1+(uv)^{N})\Big)\\
    &
    -g(1+u)^{2g-1}(1+v)^{2g-1}
    (1-(uv)^2)^2(uv)^{2g -2 +(N+1)/2}(1-(uv)^{N})
    \Bigg] \, ,
\end{aligned}
$$
where $N=d_1-d_2-2g+2$.
\end{theorem}

\begin{proof}
The first statement follows from Theorem \ref{thm:moduli(2,2)}. To
compute $e(\cN_{\s_M^+})-e(\cN_{\s_m^+})$ we use
 Proposition~\ref{prop:finally} for $\s=\s_M^+$. Note that in this case
 $n_0=[\mu_1-\mu_2]$. All the terms in the formula of
 Proposition \ref{prop:finally} involving $x^{n_0}$ yield positive
 powers of $x$, so they can be disregarded for computing
 $\coeff_{x^0}$. Hence
  $$
   \begin{aligned}
   e&(\cN_{\s_M^+}) = e(\cN_{\s_m^+}) +\coeff_{x^0}\Bigg[
    \frac{(1+u)^{3g}(1+v)^{3g}(1+ux)^{g}(1+vx)^g}
    {(1-uv)^2(1-x)(1-uvx)x^{[\mu_1-\mu_2]}} \cdot \\
    &
     \Bigg(\frac{(uv)^{g-1+(d_1-d_2-1)/2}x}{(1-(uv)^{-1}x)(1-x)}
     -\, \frac{(uv)^{(3d_1-3d_2-1)/2-g}x }{(1-(uv)^{-1}x)^2} -\, \frac{(uv)^{2g+1} x }{(1-(uv)^2x)^2}
    +\, \frac{(uv)^{d_1-d_2+2}x }{(1-(uv)^2x)(1-uvx)}\Bigg)
    \Bigg]\, .
  \end{aligned}
  $$

As $\mu_1-\mu_2>2g-2$, let $m\geq 0$ such that
$[\mu_1-\mu_2]=2g-2+m$. Introduce the following function
  $$
  F(a,b,c)=\coeff_{x^0}\Bigg( \frac{(1+ux)^{g}(1+vx)^gx^{3-2g-m}}{(1-ax)^2(1-bx)(1-cx)}
    \Bigg)=\Res_{x=0}\Bigg(\frac{(1+ux)^{g}(1+vx)^gx^{2-2g-m}}{(1-ax)^2(1-bx)(1-cx)}\Bigg)\
    ,
  $$
where $a,b,c\neq 0$. So
 \begin{equation} \label{eqn:acabando}
 \begin{aligned}
 e(\cN_{\s_M^+})=& \ e(\cN_{\smp})+
    \frac{(1+u)^{3g}(1+v)^{3g}}{(1-uv)^2} \bigg( (uv)^{3g-3+m} F(1,uv,(uv)^{-1})
    - (uv)^{5g-5+3m} F((uv)^{-1},1,uv) \\
    & - (uv)^{2g+1} F((uv)^2,1,uv)
    + (uv)^{4g-1+2m} F(uv,1,(uv)^2) \bigg)
  \end{aligned}
  \end{equation}
using $d_1-d_2=4g-3+2m$.

The function
  $$
  G(x)=\frac{(1+ux)^{g}(1+vx)^g x^{2-2g-m}}{(1-ax)^2(1-bx)(1-cx)}
  $$
is a meromorphic function on $\CC\cup \{\infty\}$ with poles at
$x=0$, $x=1/a$, $x=1/b$ and $x=1/c$. Note that there is no pole at
$\infty$. So
  $$
  F(a,b,c)=-\Res_{x=1/a} G(x) -\Res_{x=1/b} G(x) -\Res_{x=1/c}
  G(x)\, .
  $$
An easy calculation, using that
  $$
 \begin{aligned}
 &\Res_{x=1/a} G(x) = \frac{d}{dx}\bigg|_{x=1/a} \big( G(x)(x-1/a)^2\big)\,
 . \\
 &\Res_{x=1/b} G(x) = G(x)(x-1/b) |_{x=1/b}\, , \\
 &\Res_{x=1/c} G(x) = G(x)(x-1/c) |_{x=1/c}\, ,
 \end{aligned}
 $$
yields
 $$
 \begin{aligned}
    F(a,b,c)&=
    \frac{a^{m-1}b(a+u)^g(a+v)^g}{(a-b)^2(c-a)}
    +
    \frac{a^{m-1}c(a+u)^g(a+v)^g}{(b-a)(c-a)^2}\\
    &+
    \frac{b^m(b+u)^g(b+v)^g}{(a-b)^2(b-c)}
    +
    \frac{c^m(c+u)^g(c+v)^g}{(c-a)^2(c-b)}\\
    &+
    \frac{a^{m-1}(a+u)^{g-1}(a+v)^{g-1}}{(a-b)(a-c)}
    \Big(
    g \, a(2a+u+v)+(m-2)(a+u)(a+v)
    \Big) \, .
 \end{aligned}
 $$

Using this into (\ref{eqn:acabando}) and Theorem \ref{thm:(2,2,odd,even)}, we have
 $$
 \begin{aligned}
 e(\cN_{\s_M^+})
 \!=\!& \ \frac{(1+u)^{2g}(1+v)^{2g}}{(1-uv)^3(1-(uv)^2)^2}
    \Bigg[
    (1+u^2v)^{2g}(1+uv^2)^{2g}(1-(uv)^{4g+4m-2})\\
    &
    +(1-2m-2g)(1+u^2v)^g(1+uv^2)^g(1+u)^g(1+v)^g (uv)^{3g+2m-2}(1-(uv)^2)
    \\
    &
    +\!
    (1\!+\!u)^{2g}(1\!+\!v)^{2g}(1\!+\!uv)^2(uv)^{3g+m-2} 
    \Big((1\!-\!(uv)^{2g+2m}) \!-\! (m\!+\!g)(1\!-\!uv)(1\!+\!(uv)^{2g+2m-1})\Big)\\
    &
    -g(1+u)^{2g-1}(1+v)^{2g-1}
    (1-(uv)^2)^2(uv)^{3g+m-2}(1-(uv)^{2g+2m-1})
    \Bigg] \, .
\end{aligned}
$$
As $N=d_1-d_2 -2g+2 =2m+2g-1$, we get the formula in the statement.
\end{proof}

\begin{corollary}\label{cor:finally}
 Suppose that $d_1$ is even and $d_2$ is odd.
 Then the large $\s$ moduli space $\cN_{\s_M^+}=\cN^s_{\s_M^+}$ is smooth
 and compact. If $\mu_1-\mu_2>2g-2$ its Hodge polynomial has the same formula
 as that of Theorem \ref{thm:finally}.
\end{corollary}

\begin{proof}
Use the isomorphism $\cN_\s(2,2,d_1,d_2)\cong
\cN_\s(2,2,-d_2,-d_1)$ together with Theorem \ref{thm:finally}.
\end{proof}

\begin{corollary}\label{cor:finally-n}
 Suppose that $d_1+d_2$ is odd and $\mu_1-\mu_2>2g-2$. Then the
 Poincar{\'e} polynomial of $\cN_{\s_M^+}$ is
  $$
 \begin{aligned}
 P_t(\cN_{\s_M^+})=& \ \frac{(1+t)^{4g}}{(1-t^2)^3(1-t^4)^2}
    \Bigg[
    (1+t^3)^{4g}(1-t^{4N})
     - N \, (1+t^3)^{2g}(1+t)^{2g} t^{2N+2g-2} (1-t^4) \\
    & + (1+t)^{4g}(1+t^2)^2 t^{N+4g-3}
    \Big((1-t^{2N+2}) - \frac{N+1}{2} \, (1-t^2)(1+t^{2N})\Big)\\
    &
    -g \, (1+t)^{4g-2}
    (1-t^4)^2 t^{N +4g-3}(1-t^{2N})
    \Bigg]\, ,
\end{aligned}
$$
where $N=d_1-d_2 -2g+2$.
\hfill $\Box$
\end{corollary}


\begin{thebibliography}{MMMM}



\bibitem{AG} \textsc{{\'A}lvarez-C{\'o}nsul, L.; Garc{\'\i}a-Prada, O.}: Dimensional reduction,
SL (2,C)-equivariant bundles, and stable holomorphic chains.
\textsl{Internat. J. Math.} \textbf{12} (2001) 159--201.

\bibitem{AGS} \textsc{{\'A}lvarez-C{\'o}nsul, L.; Garc{\'\i}a-Prada, O.; Schmitt, A.}:
On the geometry of moduli spaces of holomorphic chains over
compact Riemann surfaces.  \textsl{Internat. Math. Res. Papers},
Art ID 73597 (2006) 1--82.

\bibitem{BD}\textsc{Bradlow, S. B.; Daskalopoulos}:
{Moduli of stable pairs for holomorphic bundles over Riemann
surfaces}. \textsl{Internat. J. Math.} \textbf{2} (1991) 477--513.

\bibitem{BGP}\textsc{Bradlow, S. B.; Garc\'{\i}a--Prada, O.}:
{Stable triples, equivariant bundles and dimensional reduction}.
\textsl{Math. Ann.} \textbf{304} (1996) 225--252.

\bibitem{BGPG}\textsc{Bradlow, S. B.; Garc\'{\i}a--Prada, O.;
Gothen, P.B}: {Moduli spaces of holomorphic triples over compact
Riemann surfaces}. \textsl{Math. Ann.} \textbf{328} (2004)
299--351.

\bibitem{Bur}\textsc{Burillo, J.}: {El polinomio de Poincar\'{e}--Hodge de un
producto sim\'{e}trico de variedades k\"{a}hlerianas compactas}.
\textsl{Collect. Math.} \textbf{41} (1990) 59--69.

\bibitem{BaR}\textsc{Del Ba{\~n}o, S.}:
On the motive of moduli spaces of rank two vector bundles over a
curve. \textsl{Compositio Math.} \textbf{131} (2002) 1--30.

\bibitem{De} \textsc{Deligne, P.:} {Th{\'e}orie de Hodge I,II,III}.
In \textsl{Proc.\ I.C.M.}, vol.\ 1, 1970, pp. 425--430; in
\textsl{Publ.\ Math.\ I.H.E.S.} \textbf{40} (1971) 5--58; ibid.\
\textbf{44} (1974) 5--77.

\bibitem{Del} \textsc{Deligne, P.}:
Th{\'e}or{\`e}me de Lefschetz et crit{\`e}res de d{\'e}g{\'e}n{\'e}rescence de suites
spectrales. \textsl{Publ.\ Math.\ I.H.E.S.} \textbf{35} (1968)
259--278.


\bibitem{Du} \textsc{Durfee, A.H.:} \textit{Algebraic
varieties which are a disjoint union of subvarieties},
\textsl{Lecture Notes in Pure Appl.\ Math.} \textbf{105}, Marcel
Dekker, 1987, pp. 99--102.

\bibitem{DK} \textsc{Danivol, V.I.; Khovanski\v{\i},
A.G.:} {Newton polyhedra and an algorithm for computing
Hodge-Deligne numbers}, \textsl{Math.\ U.S.S.R. Izvestiya}
\textbf{29} (1987) 279--298.

\bibitem{EK} \textsc{Earl, R.; Kirwan, F.:} {The Hodge numbers of the
moduli spaces of vector bundles over a Riemann surface}.
\textsl{Q. J. Math.} \textbf{51} (2000) 465--483.

\bibitem{FM}\textsc{Fulton, W.; MacPherson, R.:} {A
compactification of configuration spaces}. \textsl{Annals of
Math.} \textbf{139} (1994) 183--225.


\bibitem{GP}\textsc{Garc\'{\i}a--Prada, O.}: {Dimensional reduction of
stable bundles, vortices and stable pairs}. \textsl{Internat. J.
Math.} \textbf{5} (1994) 1--52.

\bibitem{GPGM}\textsc{Garc\'{\i}a--Prada, O.; Gothen, P.B.; Mu\~{n}oz, V.}:
{Betti numbers of the moduli space of rank $3$ parabolic Higgs
bundles}. \textsl{Memoirs Amer.\ Math.\ Soc.} In press.

\bibitem{Go} \textsc{Gothen, P.B.}: {The Betti numbers of the moduli space of stable rank
$3$ Higgs bundles on a Riemann surface}. \textsl{Internat.\ J. Math.\/}
{\bf 5} (1994) 861--875.

\bibitem{GK} \textsc{Gothen, P.B., King, A.D.}: {Homological
algebra of quiver bundles}. \textsl{J. London Math. Soc. (2)}
\textbf{71} (2005) 85--99.

\bibitem{Gri} \textsc{Griffiths, P.}: Periods of integrals on
algebraic manifolds. III. Some global differential-geometric
properties of the period mapping. \textsl{Publ.\ Math.\ I.H.E.S.}
\textbf{38} (1970) 125--180.

\bibitem{Hi} \textsc{Hitchin, N.J.}: The self-duality equations on a Riemann surface,
\textit{Proc. London Math. Soc. (3)} \textbf{55} (1987) 59--126.


\bibitem{Ki} \textsc{Kirwan, F.}: On the homology of compactifications of moduli
spaces of vector bundles over a Riemann surface.  \textit{Proc.
London Math. Soc. (3)} \textbf{53}  (1986)  237--266.


\bibitem{MOV} \textsc{Mu{\~n}oz, V.; Ortega, D.; V{\'a}zquez-Gallo, M-J.}:
Hodge polynomials of the moduli spaces of pairs.
 \textsl{Internat. J. Math.} In press.

\bibitem{NR2} \textsc{Narasimhan, M. S.; Ramanan, S.}:
Geometry of Hecke cycles. I. C. P. Ramanujam---a tribute, pp.
291--345, \textsl{Tata Inst. Fund. Res. Studies in Math.}
\textbf{8}, Springer, Berlin-New York, 1978.


\bibitem{Sch} \textsc{Schmitt, A.}:
{A universal construction for the moduli spaces of decorated
vector bundles}. \textsl{Transform. Groups} \textbf{9} (2004)
167--209.

\bibitem{Th}\textsc{Thaddeus, M.}:
{Stable pairs, linear systems and the Verlinde formula}.
\textsl{Invent. Math.} \textbf{117} (1994) 317--353.

\end{thebibliography}
\end{document}